\documentclass[a4paper]{article}

\usepackage{amsfonts}
\usepackage{amsthm}
\usepackage{amssymb}
\usepackage{amsmath}
\usepackage{upref}
\usepackage{mathrsfs}
\usepackage{color}
\usepackage[citecolor=blue,colorlinks=true]{hyperref}
\usepackage{enumitem}

\newtheorem{theorem}{Theorem}
\newtheorem{proposition}[theorem]{Proposition}
\newtheorem{lemma}[theorem]{Lemma}

\newtheorem{definition}[theorem]{Definition}
\newtheorem{remark}[theorem]{Remark}

\def\N{\mathbb{N}}
\def\Z{\mathbb{Z}}

\def\R{\mathbb{R}}
\def\C{\mathbb{C}}
\def\S{\mathbb{S}}

\def\carac{\mathbf{1}}

\def\ds{\displaystyle} 
\def\div{{\rm div}}
\def\refe#1{(\ref{#1})}

\def\ocirc#1{\ifmmode\setbox0=\hbox{$#1$}\dimen0=\ht0
    \advance\dimen0 by1pt\rlap{\hbox to\wd0{\hss\raise\dimen0
    \hbox{\hskip.2em$\scriptscriptstyle\circ$}\hss}}#1\else
    {\accent"17 #1}\fi} 
    
\def\qed{\rule{0.2cm}{0.2cm}}

\def\eps{\varepsilon}

\def\<{\langle}
\def\>{\rangle}
\def\F{\mathcal{F}}
\def\GG{\mathbf{G}}
\def\P{\mathbb{P}}
\def\E{\mathbb{E}}
\def\T{\mathbb{T}}

\DeclareMathOperator{\esssup}{ess\,sup}

\begin{document}

\title{Invariant measure of scalar first-order conservation laws with stochastic forcing}
\author{A. Debussche\thanks{IRMAR, ENS Cachan Bretagne, CNRS, UEB. av Robert Schuman, F-35170 Bruz, France. Email: arnaud.debussche@bretagne.ens-cachan.fr, partially supported by ANR STOSYMAP} and J. Vovelle\thanks{Universit\'e de Lyon ; CNRS ; Universit\'e Lyon 1, Institut Camille Jordan,  43 boulevard du 11 novembre 1918, F-69622 Villeurbanne Cedex, France. Email: vovelle@math.univ-lyon1.fr, partially supported by ANR STOSYMAP and ANR STAB}}
\maketitle

\begin{abstract} Under an hypothesis of non-degeneracy of the flux, we study the long-time behaviour of periodic scalar first-order conservation laws with stochastic forcing in any space dimension. 
For sub-cubic fluxes, we show the existence of an invariant measure. Moreover for sub-quadratic 
fluxes we show uniqueness and ergodicity of the invariant measure. Also, since this invariant 
measure is supported by
$L^p$ for some $p$ small, we are led to generalize to the stochastic case the theory of $L^1$ solutions developed by Chen and Perthame \cite{ChenPerthame03}.
\end{abstract}

\medskip


\tableofcontents

\section{Introduction}


In this work we investigate the long-time behaviour of periodic scalar first-order conservation laws with stochastic forcing. In space dimension one, there is a famous paper of E, Khanin, Mazel, Sinai \cite{EKMS00} where the author prove existence and uniqueness, and also analyse the invariant measure for the periodic inviscid Burgers equation with stochastic forcing. The analysis is done by use of the Lax-Oleinik formula, which means that the random unknown $u$ is derived from a potential $\psi$ which satisfies a periodic Hamilton-Jacobi with stochastic forcing. This can be extended to higher dimension (Iturriaga and Khanin \cite{IturriagaKhanin03}) and has also been extended to the case of fractional noise by Saussereau and Stoica, \cite{SaussereauStoica12}. Recently, Boritchev (\cite{boritchev}) has been able 
to derive sharp estimates on the solutions of the stochastic Burgers equation in the viscous case. They hold independently of the viscosity (see Remark~\ref{rkBoritchev} on that point). 
\smallskip

Note that in all these articles, the space domain is compact. A recent work by Bakhtin \cite{Bakhtin13} deals with a scalar first-order conservation law with Poisson random forcing set on the whole line. We also would like to mention, for
 stochastically forced Hamilton-Jacobi equations, the thorough analysis of the invariant measure for such problems by Dirr and Souganidis in \cite{DirrSouganidis05}.

\medskip

Our purpose is to generalize \cite{EKMS00} to higher dimension and, mainly, to relax the hypothesis of uniform convexity which is necessary when the Lax-Oleinik formula is used. 
The first step was to have a satisfactory framework for existence and uniqueness of solutions. This has 
been accomplished by several authors and this point is now rather understood (see Kim \cite{kim}, Vallet Wittbold \cite{ValletWittbold09}, Feng, Nualart, \cite{FengNualart08},  Debussche, Vovelle \cite{DebusscheVovelle10,DebusscheVovelle10revised}, Bauzet, Vallet, Wittbold \cite{bauzet}, Chen, Ding, Karlsen \cite{ChenDingKarlsen12}, Lions, Perthame, Souganidis \cite{LionsPerthameSouganidis13}). 
\medskip

In this article we aim at proving existence and uniqueness of invariant measures for the type of stochastic conservation laws studied in \cite{DebusscheVovelle10}. There, a general result of existence and uniqueness result was obtained thanks to the kinetic formulation (see Lions, Perthame, Tadmor \cite{lions}, Perthame \cite{perth}) suitably generalized to the stochastic case. 
This kinetic formulation seems well adapted for this purpose since it allows to keep track of the dissipation due to the shocks. Note that, as explained below, we need 
here a notion of solution in $L^1$, which we develop in this article in the framework of kinetic formulation.
\medskip

Under an hypothesis of non-degeneracy of the flux, we show  the existence of an invariant measure in any dimension of space. For technical reasons we need to assume that the flux 
does not grow faster than a cubic polynomial. Moreover,  for sub-quadratic fluxes we show the uniqueness of the invariant measure, see our main result, Theorem~\ref{th:maintheorem}. An essential 
tool for the proofs is the use of averaging lemma for kinetic equations. We use such an averaging 
lemma well adapted to stochastic equation issued from Bouchut and Desvillettes \cite{BouchutDesvillettes99}.

\medskip

More precisely, let $(\Omega,\F,\P,(\F_t),(\beta_k(t)))$ be a stochastic basis and let $T>0$. We study the invariant measure for the first-order scalar conservation law with stochastic forcing
\begin{equation}\label{stoSCL}
du+\div(A(u))dt=\Phi dW(t),\quad x\in\T^N, t\in(0,T).
\end{equation}
The equation is periodic in the space variable $x$:  $x\in\T^N$ where $\T^N$ is the $N$-dimensional torus. The flux function $A$ in \refe{stoSCL} is supposed to be of class $C^2$: $A\in C^2(\R;\R^N)$ and its derivatives have at most polynomial growth. We assume that the filtration $(\mathcal F_t)_{t\ge 0}$ is complete and that $W$ is a cylindrical Wiener process: $W=\sum_{k\geq 1}\beta_k e_k$, where the $\beta_k$ are independent Brownian processes and $(e_k)_{k\geq 1}$ is a complete orthonormal system in a Hilbert space $H$. The map $\Phi\colon H\to L^2(\T^N)$ is defined by $\Phi e_k=g_k$ where $g_k$ is a regular function on $\T^N$. More precisely, we assume
$g_k\in C(\T^N)$, with the bounds 
\begin{align}
\GG^2(x)=\sum_{k\geq 1}|g_k(x)|^2\leq D_0,\label{D0}\\
\sum_{k\geq 1}|g_k(x)-g_k(y)|^2\leq D_1|x-y|^2,\label{D1}
\end{align}
for all $x,y\in\T^N$. Note in particular that $\Phi\colon H\to L^2(\T^N)$ is Hilbert-Schmidt since $\|g_k\|_{L^2(\T^N)}\leq\|g_k\|_{C(\T^N)}$ and thus
\begin{equation*}
\sum_{k\geq 1}\|g_k\|_{L^2(\T^N)}^2\leq D_0.
\end{equation*}
Existence and uniqueness of a solution $u^\omega(t)$ to \refe{stoSCL} satisfying a given initial 
condition $u(0)=u_0\in L^\infty(\T^N)$ has been proved in \cite{DebusscheVovelle10}. This result can 
easily be extended to initial data in $L^p(\T^N)$ for $p$ larger than the degree of polynomial growth of 
$A$. However, as stated below, the invariant measure is supported by 
$L^r(\T^N)$ with $r$ small and we do not know whether $r$ can be taken sufficiently large so that 
the result of \cite{DebusscheVovelle10} cannot be used for such invariant measures. Thus, we need to 
develop a theory of existence and uniqueness for initial data in $L^r(\T^N)$ for small $r$. In fact, we 
generalize to the stochastic context the notion of solutions in $L^1(\T^N)$ developed in 
\cite{ChenPerthame03} . Note however that in \cite{ChenPerthame03} second-order, possibly degenerate equations are considered, so that, strictly speaking, what we generalize is the ``fully degenerate" case contained in \cite{ChenPerthame03}. Anyway, what matter technically here is to prove that the kinetic measure
decay sufficiently fast at infinity. 
\medskip

Let us now give some precisions on our framework. We assume that the noise is additive and that the functions $g_k$ satisfy the cancellation condition
\begin{equation}
\forall k\geq 1,\quad \int_{\T^N} g_k(x)dx=0.
\label{intgk}\end{equation}
The solution then satisfies $\bar u(t)=\int_{\T^N}u(t,x) dx = \int_{\T^N} u(0,x)dx$ almost surely. We will consider initial data with non random space average $\bar u$. Then this remains the case for all time 
$t\ge 0$. Changing $u$ to $u-\bar u$, we see that it is no loss of generality to consider the 
case $\bar u=0$. Thus in all the article, we only consider such initial data.
\medskip

To introduce the non-degeneracy condition on which we will work on, let us set
\begin{equation}\label{defiota}
\iota(\eps)=\sup_{\alpha\in\R,\beta\in\S^{N-1}} |\{\xi\in\R; |\alpha+\beta\cdot a(\xi)|<\eps\}|,
\end{equation}
and
\begin{equation}\label{defeta}
\eta(\eps):=\int_0^\infty e^{-t}\iota(t\eps) dt.
\end{equation}
It is known that, if $a:=A'$ satisfies the following non-degeneracy condition:
\begin{equation}
\lim_{\eps\to 0}\iota(\eps)=0,
\label{and1}\end{equation}
then all solutions of the deterministic equation, {\it i.e.} with $\Phi=0$, converge to zero (\textit{cf.} Debussche, Vovelle \cite{DebusscheVovelle09} for example, in that case the condition can be localized. See also
\cite{chen-perthame09} for the quasilinear parabolic case). The condition~\eqref{and1} is a condition of non-stationarity of $\xi\mapsto a(\xi)$. 
\medskip


In this paper, we use the approach developed in \cite{BouchutDesvillettes99}. There, an averaging 
Lemma based on the Fourier transform in $x$ is developed. The non-degeneracy assumptions is 
strengthened into:
\begin{equation}
\label{andiota}
\iota(\eps)\le c_1 \eps^b,
\end{equation}
for some $c_1>0$ and $b>0$. Note that $b\leq 1$, unless $a\equiv 0$, and that \refe{andiota} is equivalent  to 
\begin{equation}
\label{and}
\eta(\eps)\le c_1 \eps^b,
\end{equation}
for possibly a different constant $c_1$. We will use the averaging lemma (under the form developed in \cite{BouchutDesvillettes99}, \textit{cf.} Section~\ref{sec:bound} and Section~\ref{sec:uniqueness}) to prove the following result.

\begin{theorem}[Invariant measure] Let $A\in C^2(\R;\R^N)$ satisfy the non-de\-ge\-ne\-ra\-cy condition \refe{andiota} where $\iota$ is defined by \refe{defiota}. Assume conditions \refe{D0}-\refe{D1} on the noise, the cancellation condition \refe{intgk} and that $A$ is at most cubic in the following sense:
\begin{equation}\label{aquad}
|a'(\xi)|\le c_1(|\xi|+1),\; \xi\in \R.
\end{equation}
Then {\rm there exists an invariant measure} for \refe{stoSCL} in $L^1(\T^N)$.  Moreover, it is supported
by $L^r(\T^N)$ for $r<2+\frac{b}2$ if $N=1$ and $r<\frac{N}{N-1}$ if $N\ge 2$. If the condition \refe{aquad} is strengthened into the hypothesis that $A$ is sub-quadratic in the following sense:
\begin{equation}\label{linear}
|a'(\xi)|\le c_2,\quad \xi \in \R,
\end{equation}
then {\rm the invariant measure is unique}.
\label{th:maintheorem}\end{theorem}

Thanks to the kinetic formulation, it is easy to see that  a stationary solution develops shocks (see Remark \ref{rmk8}). Thus the noise does not have a regularizing effect here (see Flandoli \cite{flandoli11} for a  review of such effects in SPDEs). Note that the presence of
shocks was already observed in \cite{EKMS00} and their structure was carefully studied.
\medskip

In the statement of Theorem~\ref{th:maintheorem} and in the sequel, we use the letter $C$ or $c$ to denote a constant. Its value may change from one line
to another. Sometimes, we precise its dependence on some parameters. 
\medskip

\begin{remark} The existence of an invariant measure is closely related to existence of uniform in time estimates for the solution to \eqref{stoSCL} (we do not specify in which norm in this discussion). Note that such uniform estimates, with respect to time and with respect to the viscosity parameter also since second-order equations are considered, have been given by Boritchev, \cite{boritchev}, for the generalized Burgers equation with a flux $A$ satisfying an hypothesis of strict convexity and an hypothesis of (strict) sub-quadratic growth of $A'$. Compare to \eqref{aquad} here.
\label{rkBoritchev}\end{remark}

To prove Theorem~\ref{th:maintheorem}, we first describe in Section~\ref{sec:kisol}  the concept of $L^1$ solution for \refe{stoSCL} used in this article. The well posedness in $L^1$ is proved in Appendix~\ref{app:L1}. Then, in Section~\ref{sec:bound}, we prove some bounds uniform in time on the solution which give the existence of the invariant measure. These bounds are used again in Section~\ref{sec:uniqueness}, together with a quite classical argument of smallness of the noise, to obtain the uniqueness of the invariant measure. The difficulty in the proof 
of uniqueness is that, contrary to the convex case treated in the references above, we do not have
any uniformity with respect to the initial data on the proximity between the deterministic solution and a solution with small noise. A similar thing appears for the 2D stochastic Navier-Stokes equations with large viscosity (see Mattingly \cite{mattingly99}) but in this case the a priori estimates are easier to obtain and the contraction 
of the trajectories is much stronger.


\section{Kinetic solution}\label{sec:kisol}



Let us give the notion of solution  to \refe{stoSCL} we need use in this article.  It has been introduced in \cite{DebusscheVovelle10}. We slightly modify the definition to adapt the $L^1$ setting.

\begin{definition}[Kinetic measure] We say that a map $m$ from $\Omega$ to the set of non-negative Radon measures over $\T^N\times[0,T]\times\R$ is a kinetic measure if
\begin{enumerate}
\item $m$ is measurable, in the sense that for each $\phi\in C_c(\T^N\times[0,T]\times\R)$, $\<m,\phi\>\colon\Omega\to\R$ is,
\item $m$ vanishes for large $\xi$ in the sense that
\begin{equation}
\lim_{n\to+\infty}\E m(A_n)=0,
\label{inftymL1}\end{equation}
where $A_n=\T^N\times[0,T]\times\{\xi\in\R,n\leq|\xi|\leq n+1\}$,
\item for all $\phi\in C_c(\T^N\times\R)$, the process
\begin{equation*}
t\mapsto \int_{\T^N\times[0,t]\times\R} \phi(x,\xi)dm(x,s,\xi)
\end{equation*}
is predictable.
\end{enumerate}
\label{def:kineticmeasureL1}\end{definition}

\begin{definition}[Solution] Let $u_{0}\in L^1(\T^N)$. A measurable function $u\colon\T^N\times [0,T]\times\Omega\to\R$ is said to be a solution to~\refe{stoSCL} with initial datum $u_0$ if $(u(t))$ is predictable, if 
\begin{equation}
\E\left(\esssup_{t\in[0,T]}\|u(t)\|_{L^1(\T^N)}\right)<+\infty,
\label{eq:integrabilityuL1}\end{equation} 
and if there exists a kinetic measure $m$ such that $f:=\mathbf{1}_{u>\xi}$ satisfies:  for all $\varphi\in C^1_c(\T^N\times[0,T)\times\R)$, 
\begin{multline}
\int_0^T\<f(t),\partial_t \varphi(t)\>dt+\<f_0,\varphi(0)\>
+\int_0^T \<f(t),a(\xi)\cdot\nabla\varphi(t)\>dt\\
=-\sum_{k\geq 1}\int_0^T\int_{\T^N}g_k(x)\varphi(x,t,u(x,t)) dxd\beta_k(t)\\
-\frac{1}{2}\int_0^T\int_{\T^N} \partial_\xi\varphi(x,t,u(x,t))\GG^2(x) dx dt+m(\partial_\xi\varphi).
\label{eq:kineticupreL1}\end{multline}
\label{defkineticsolL1}\end{definition}

In \refe{eq:kineticupreL1}, $f_0(x,\xi)=\mathbf{1}_{u_0(x)>\xi}$. We have used the brackets $\<\cdot,\cdot\>$ to denote the duality between $C^\infty_c(\T^N\times\R)$ and the space of distributions over $\T^N\times\R$. In what follows, we will denote similarly the integral 
\begin{equation*}
\<F,G\>=\int_{\T^N}\int_\R F(x,\xi)G(x,\xi) dx d\xi,\quad F\in L^p(\T^N\times\R), G\in L^q(\T^N\times\R),
\end{equation*}
where $1\leq p\leq +\infty$ and $q$ is the conjugate exponent of $p$. In \refe{eq:kineticupreL1} also, we have used the shorthand $m(\phi)$ for
\begin{equation*}
m(\phi)=\int_{\T^N\times[0,T]\times\R} \phi(x,t,\xi)dm(x,t,\xi),\quad \phi\in C_b(\T^N\times[0,T]\times\R).
\end{equation*}

Switching from $\mathbf{1}_{u>\xi}$ to $\chi_u(\xi):=\mathbf{1}_{u>\xi>0}-\mathbf{1}_{0>\xi>u}=\mathbf{1}_{u>\xi}-\mathbf{1}_{0>\xi}$, equation~\refe{eq:kineticupreL1} is the weak form of the equation
\begin{equation}
(\partial_t+a(\xi)\cdot\nabla)f=\delta_{u=\xi}\dot W+\partial_\xi(m-\frac{1}{2}\GG^2\delta_{u=\xi}),\quad f=\chi_u.
\label{eq:kineticD}\end{equation}

In Appendix~\ref{app:L1} we prove the following result which generalizes  \cite{DebusscheVovelle10} ({\it cf.} Theorem~11, Corollary~12, and Theorem~19 in \cite{DebusscheVovelle10} in particular for a precise definition of the ``parabolic approximation").

\begin{theorem}[Resolution of \refe{stoSCL} in $L^1$] Let $u_0\in L^1(\T^N)$. There exists a unique measurable $u\colon\T^N\times [0,T]\times\Omega\to\R$ solution to~\refe{stoSCL} with initial datum $u_0$ in the sense of Definition~\ref{defkineticsolL1}. Besides, $u$ has almost surely continuous trajectories in $L^1(\T^N)$ and $u$ is the a.s. limit in $L^1(\T^N\times(0,T))$ of the parabolic approximation to \refe{stoSCL}.  Moreover, given $u_0^1$ and $u_0^2\in L^1(\T^N)$, the following holds:
\begin{equation}\label{contractionL1}
\|u^1(t)-u^2(t)\|_{L^1(\T^N)}\le \|u_0^1-u_0^2\|_{L^1(\T^N)},\; a.s.
\end{equation}
\label{th:ResolutionL1}\end{theorem}

We need this generalization because we are not able to prove that the invariant measure we construct
has its support in $L^p(\T^N)$ for $p$ sufficiently large to apply the result in 
 \cite{DebusscheVovelle10}. Note that it follows from the proof that if $u_0\in L^r(\T^N)$ for some $r>1$ then the solution lives in $L^r(\T^N)$ and an estimate similar to \eqref{eq:integrabilityuL1} but in $L^r(\T^N)$ holds. 
 
 \medskip
By Theorem~\ref{th:ResolutionL1}, we can define the transition semigroup in $L^1(\T^N)$:
$$
P_t\phi(u_0)=\E(\phi(u(t))),\; \phi\in \mathcal B_b(L^1(\T^N)).
$$
It follows easily from the property of $L^1$-contraction \eqref{contractionL1} that $(P_t)$ is Feller.
\medskip

Our aim is to construct an invariant measure for $(P_t)$ under assumption \eqref{aquad}. We prove the first part of Theorem~\ref{th:maintheorem} (existence of an invariant measure) in the next section and we conclude this section with some remarks.
\bigskip

Since the solution of \eqref{stoSCL} has continuous trajectories, it is easy to prove that
it satisfies the following weak formulation:
\begin{multline}
-\<f(t),\varphi\>dt+\<f_0,\varphi\>
+\int_0^t \<f(s),a(\xi)\cdot\nabla\varphi\>ds\\
=-\sum_{k\geq 1}\int_0^t\int_{\T^N}g_k(x)\varphi(x,u(x,t)) dxd\beta_k(s)\\
-\frac{1}{2}\int_0^t\int_{\T^N} \partial_\xi\varphi(x,u(x,s))\GG^2(x) dx dt+m(\partial_\xi\varphi),
\label{eq:kinetic-weak}\end{multline}
 for all $\varphi\in C^1_c(\T^N\times\R)$.
\medskip
 
In fact, we will even use stronger formulation by using the following remark.

\begin{remark}[Regularity of the parabolic approximation] Let the pa\-ra\-bo\-lic approximation to \eqref{stoSCL} be given by:
\begin{equation}\label{parabolic}
du^\eta+\div(A_\eta(u))dt-\eta \Delta u^\eta=\Phi_\eta dW(t),\quad x\in\T^N, t\in(0,T),
\end{equation}
where $A_\eta$ and $\Phi_\eta$ are smooth. By Hofmanov\'a \cite{Hofmanova13}, if $u^\eta(0)$
is in $W^{m,p}(\T^N)\cap W^{1,mp}(\T^N)$, then $u^\eta$ belongs to $L^q(\Omega;C([0,T];W^{m,p}(\T^N)\cap W^{1,mp}(\T^N)))$ for any $q\ge 1$. Also by Theorem~\ref{th:ResolutionL1}, if one approximates the solution of \eqref{stoSCL} by the solution of \eqref{parabolic} with $u^\eta(0) \to u_0$ in $L^1(\T^N)$, then $u^\eta\to u$ in $L^\infty(0,T;L^1(\Omega\times\T^N))$.
\label{rk:regParabApprox}\end{remark}

\begin{remark}[Multiplicative noise] A natural extension of our study is  to  consider a noise depending on the solution of the form 
$$
\Phi(u)dW=\sum_{k\ge 1} \Phi(u)e_kd\beta_k=\sum_{k\ge 1} g_k(u)d\beta_k.
$$
With such noise, the spatial average satisfies:
$$
\int_{\T^N} u(x,t)dx = \int_{\T^N} u_0(x)dx +\sum_{k\ge 1}\int_0^t \int_{\T^N}  g_k(u(x,s)dx d\beta_k(s).
$$
It seems difficult to provide any bound on this average for a general noise. A first possibility is to consider a noise satisfying
$$
\int_{\T^N}g_k(u(x))dx=0,\; k\ge 1.
$$
Then the mass is preserved exactly. This would hold for $g_k(u)=\frac{\partial}{\partial x}h_k(u)$ or
$g_k(u)=h_k(u)-\frac1{|\T^N|}\int_{\T^N} h_k(u)dx.$ But such noises are not suited to the
framework of kinetic solutions and are not covered by our theory. 
The first one is considered in \cite{LionsPerthameSouganidis13} but the longtime behaviour is not analysed there. 
The second noise is non local and does not seem to be justified physically. 
\smallskip

It would be also possible to consider noises of the form $u(1-u)dW$. Then, the solution is easily
proved to live in $[0,1]$ and the existence of an invariant measure is easy. The longtime behaviour
in this case is expected to be completely different (see Ikeda, Watanabe \cite{ikeda-watanabe} for the case of an SDE and
Berg\'e, Saussereau \cite{berge-saussereau} for the case of a parabolic SPDE). 
\end{remark}
\begin{remark}[Shocks do not disappear]
Once we have proved Theorem \ref{th:maintheorem}, it is not difficult to construct a stationary solution of
\refe{stoSCL} with the invariant law given by this result. It is of course a kinetic solution. Assume that $N=1$, then the invariant measure is supported by $L^p(\T^N)$ for some $p>2$. It is easy to see that a preliminary step of truncation allows to take the test function $\varphi(\xi)=\xi$ in \refe{eq:kinetic-weak} and we obtain:
$$
\E(\|u(t+1)\|_{L^2(\T^N)}^2)+\E m(\T^N\times\R\times [t_0,t_{0+1}])=\E(\|u(t)\|_{L^2(\T^N)}^2)
+D_0.
$$
Thus, by stationarity,
$$
\E m(\T^N\times\R\times [t_0,t_{0+1}])=D_0.
$$
As expected, the lhs does not depend on $t_0$. More interestingly, it is not $0$. This shows that
shocks are present in the stationary solutions and noise has no regularizing effect. If $N\ge 2$, one can use the test functions of the proof of Proposition \ref{prop:decaym} to see that again shocks are present for a stationary solution. 
\label{rmk8}
\end{remark}
\section{Uniform bound and tightness for the stochastically forced equation}\label{sec:bound}

\subsection{Decomposition of $u$}\label{sec:decu}

Let $\alpha, \gamma,\delta>0$ be some positive parameters. On $L^2(\T^N)$, we consider the following regularization of the operator $-a(\xi)\cdot\nabla$:
$$
A_\gamma:=-a(\xi)\cdot\nabla-B_\gamma,\quad A_{\gamma,\delta}=A_\gamma -\delta\mathrm{Id},\quad B_\gamma:=\gamma(-\Delta)^{\alpha}.
$$
Let also $S_{A_\gamma}(t)$ and $S_{A_{\gamma,\delta}}(t)$ be the associated semi-groups on $L^2(\T^N)$:
$$
S_{A_\gamma}(t)v(x)=(e^{-tB_\gamma}v)(x-a(\xi)t),\quad S_{A_{\gamma,\delta}}(t)v(x)=e^{-\delta t}(e^{-tB_\gamma}v)(x-a(\xi)t).
$$
Note that the semi-group $e^{-tB_\gamma}$ has the following regularizing properties, which we will use later.
\begin{lemma} For $\gamma>0$, $\alpha\in(0,1]$, let $B_\gamma=\gamma(-\Delta)^{\alpha}$ on $\T^N$ and let $e^{-tB_\gamma}$ denote the associated semi-group. Then there exists a constant $c_3$ depending on $N,\, n,\,  m,\,  \alpha,\,  \beta$ such that
$$
\|(-\Delta)^{\beta/2}e^{-tB_\gamma}\|_{L^m\to L^n}\leq \frac{c_3}{(\gamma t)^{\frac{N}{2\alpha}\left(\frac1m-\frac1n\right)+\frac\beta{2\alpha}}},
$$
for all $1\leq m\leq n\leq +\infty$ and $\beta\geq 0$.
\label{lemmahalfDelta}\end{lemma}

{\bf Proof:} it is sufficient to consider the case $\gamma=1$. Note that $e^{-tB_1}\varphi=K_{t}*\varphi$ for all $\varphi\in\mathcal{C}^2(\T^N)$, where
$$
K_t(x)=\sum_{n\in\Z^d}e^{-t|n|^{2\alpha}+2\pi i n\cdot x}.
$$
Let $\mathcal{F}$ denote the Fourier transform 
$$
\mathcal{F}u(\xi)=\int_{\R^N}u(x)e^{-2\pi ix\cdot \xi}dx,
$$
for $u\in\mathcal{S}(\R^N)$. Let $H_t$ denote the inverse Fourier transform of $\xi\mapsto e^{-t|\xi|^{2\alpha}}$ (this is the ``heat" kernel associated to $(-\Delta)^{\alpha}$ on $\R^N$) and let $H_t^\mathrm{per}$ denote the periodic function generated by $H_t$:
$$
H_t^\mathrm{per}=\sum_{l\in\Z^N}H_t(x+l),\quad x\in\R^N.
$$
For a fixed $x\in\R^N$, we apply the Poisson formula
$$
\frac{1}{(2\pi)^N}\sum_{k\in\Z^N}\mathcal{F}\theta(k)=\sum_{l\in\Z^N}\theta(2\pi l), \quad\theta\in\mathcal{S}'(\R^N)
$$
to the function $\theta(z)=H_t(x+z/2\pi)$. We obtain $H_t^\mathrm{per}=K_t$.
In particular, we have
$$
\|K_t\|_{L^r(\T^N)}\leq \|H_t\|_{L^r(\R^N)},\quad 1\leq r\leq +\infty.
$$ 
By the Young inequality, we also have
$$
\|e^{-tB_1}\|_{L^m(\T^N)\to L^n(\T^N)}\leq \|K_t\|_{L^r(\T^N)},\quad \frac1m+\frac1r=1+\frac1n,
$$
and, therefore,
$$
\|e^{-tB_1}\|_{L^m(\T^N)\to L^n(\T^N)}\leq \|H_t\|_{L^r(\R^N)}.
$$
By homogeneity, $\|H_t\|_{L^r(\R^N)}=\frac{c}{t^{\frac{N}{2\alpha}\left(\frac1m-\frac1n\right)}}$. The case $\beta\not=0$ is similar. \qed
\bigskip

Now, we formally decompose the solutions using the following rewriting of  \refe{eq:kineticD}:
 \begin{equation}
(\partial_t-A_{\gamma,\delta})f=(B_\gamma+\delta \mbox{Id}) + p+\partial_\xi q,\quad f=\chi_u.
\label{eq:kineticD-1}\end{equation}
with $p= \delta_{u=\xi}\dot W$, $q=m-\frac{1}{2}\GG^2\delta_{u=\xi}$.
\medskip

More precisely, by a preliminary step of regularization, we may use the test function
$S^*_{A_{\gamma,\delta}}(T-t)\varphi$ with $\varphi \in C(\T^N)$ in \eqref{eq:kineticupreL1}.
Then, by the commutation identity
$$
S_{A_{\gamma,\delta}}(t)\partial_\xi g=\partial_\xi \left(S_{A_{\gamma,\delta}}(t)g\right)+tXS_{A_{\gamma,\delta}}(t) g,\quad X=a'(\xi)\cdot\nabla,
$$ 
 we obtain the following decomposition of the solution:
\begin{equation}\label{decompusto}
u=u^0+u^\flat+P+Q,
\end{equation}
where
\begin{align}
u^0(t)&=\int_\R S_{A_{\gamma,\delta}}(t)f(0,\xi)d\xi,\label{defu00}\\
u^\flat(t)&=\int_\R \int_0^t S_{A_{\gamma,\delta}}(s)(B_\gamma f+\delta f)(t-s,\xi)ds d\xi,\label{defuflat}\\
\<P(t),\varphi\>&=\sum_{k\geq 1} \int_{\T^N}\int_0^t \left(S_{A_{\gamma,\delta}}^*(t-s)\varphi\right)(x,u(s,x)) g_k(x) d\beta_k(s)dx,\label{defPsto}\\
\mathrm{and}\nonumber\\
\<Q(t),\varphi\>&=\int_{\T^N\times[0,t]\times\R} (t-s)a'(\xi)\cdot\nabla S_{A_{\gamma,\delta}}^*(t-s)\varphi dq(x,s,\xi),\label{defQsto}
\end{align}
where $\varphi\in C(\T^N)$ and $\<M,\varphi\>$ denote the duality product between the space of finite Borel measures on $\T^N$ and $C(\T^N)$.
\medskip

We estimate each term separately. In fact, we estimate the parabolic approximation of $u$. 
Since it is as smooth as we need (\textit{cf.} Remark~\ref{rk:regParabApprox}), the computations below are easily justified.  In particular, for 
the parabolic approximation, the kinetic measure is given by 
$$
m^\eta=\eta |\nabla u^\eta|^2 \delta_{u^\eta=\xi}
$$
which is smooth in $t$ and $x$. To lighten the computation below, we will however omit to write the dependence on $\eta$, except in section 
\ref{conclusion} where we derive the final estimate on the true solution $u$.

\subsection{Estimate of $u^0$}\label{u0}

We use the Fourier transform with respect to $x\in\T^d$, \textit{i.e.}
$$
\hat v(n)=\int_{\T^N} v(x) e^{-2\pi i n\cdot x},\quad n\in \Z^N,\quad v\in L^1(\T^N).
$$ 
After Fourier transform, $S_{A_{\gamma,\delta}}(t)$ is multiplication by $e^{-(ia(\xi)\cdot n +\gamma |n|^{2\alpha}+\delta)t}$, hence, for all $n\in \Z^N$, $n\ne 0$,
$$
\hat u^0(n,t) =\int_\R e^{-(ia(\xi)\cdot n +\gamma |n|^{2\alpha}+\delta)t}\widehat{f_0}(n,\xi)d\xi.
$$
Let us set, for $\varphi\,:\, \Z^N\times \R\to \C$ satisfying $\varphi(n,\cdot)\in L^2(\R)$ for all
$n\in \Z^N$:
\begin{equation}\label{defGG}
{\mathcal G}\varphi(n,z)=\int_\R e^{-ia(\xi)\cdot z} \varphi(n,\xi)d\xi.
\end{equation}
The operator $\mathcal{G}$ is well-defined by \cite{BouchutDesvillettes99}. Let us also define 
\begin{equation}\label{defomegan}
\omega_n=\gamma |n|^{2\alpha -1} +\delta |n|^{-1}.
\end{equation}
Then, for $T\ge0$, $n\ne0$,
\begin{align*}
\int_0^T|\hat u^0(n,t)|^2 dt &\ds = \int_0^T e^{-2(\gamma |n|^{\alpha}+\delta )t} \left| {\mathcal G}\widehat f_0 (n,nt)\right|^2 dt\\
& \le \frac1{|n|}\int_{\R^+} e^{-2\omega_n s} \left| {\mathcal G}\widehat f_0 \left(n,\frac{n}{|n|}s\right)\right|^2 ds,
\end{align*}
thanks to the change of variable $s=|n| t$.  We now use Lemma 2.4 in \cite{BouchutDesvillettes99}, which gives an estimate of the oscillatory integral \eqref{defGG}, and the obvious inequality $e^{-a}\le 1/(1+a^2)$ to deduce:
\begin{align*}
\int_0^T|\hat u^0(n,t)|^2 dt &\ds \le  \frac1{|n|}\int_{\R^+}  \frac{1}{1+4\omega_n^2s^2}  \left| {\mathcal G}\widehat f_0 \left(n,\frac{n}{|n|}s\right)\right|^2 ds\\
&\le  \frac1{|n|}\frac{2\pi}{2\omega_n}\eta(\omega_n) \left\| \widehat f_0(n,\cdot)\right\|_{L^2_\xi}^2.
\end{align*}
Recall that $\eta$ was defined in \eqref{defeta} and satisfies the decay condition \eqref{and}. Consequently, we have the estimate
$$
 \int_0^T|\hat u^0(n,t)|^2 dt \le c\frac{1}{|n|}\omega_n^{b-1}\left\| \widehat f_0(n,\cdot)\right\|_{L^2_\xi}^2.
$$
Clearly,
$$
|n|\omega_n^{1-b}\ge \gamma^{1-b} |n|^{2\alpha(1-b)+b},
$$
for $n\not=0$. Summing over $n\in\Z^N$ and noticing that 
\begin{align*}
\hat u^0(0,t)= \int_{\T^N} u^0(t,x) dx& = \int_\R \int_{\T^d} e^{-\delta t} f(0,x,\xi) d\xi dx\\
& = \int_{\T^d} e^{-\delta t} u_0(x) dx=0,
\end{align*} we obtain:
$$
\int_0^T\|u^0(t)\|_{H^{\alpha+\left(\frac12-\alpha\right)b}}^2 dt \le c\gamma^{b-1}\left\| f_0\right\|_{L^2_{x,\xi}}^2.
$$
Since $\left\| f_0\right\|_{L^2_{x,\xi}}^2= \left\| f_0\right\|_{L^1_{x,\xi}}= \left\| u_0\right\|_{L^1_{x}}$, this gives
\begin{equation}\label{estim-u0}
 \int_0^T\|u^0(t)\|_{H^{\alpha+\left(\frac12-\alpha\right)b}}^2 dt \le c\gamma^{b-1}\left\| u_0\right\|_{L^1_{x}}.
\end{equation}

\subsection{Estimate of $u^\flat$}\label{uflat}

Recall that 
$$
u^\flat(t)=\int_\R \int_0^t S_{A_{\gamma,\delta}}(s)(B_\gamma f+\delta f)(t-s,\xi)ds d\xi.
$$

We use the same argument as above. We first notice that 
$$
\widehat{u}^\flat(0,t)=\int_{\T^N} u^\flat(x,t)dx=0
$$ 
for all
$t\ge 0$, almost surely. This is easily seen using that
$$
\int_{\T^N}\int_\R f(x,t,\xi)d\xi \, dx=0,
$$ 
which 
follows from \eqref{intgk} and \eqref{eq:kinetic-weak} with test functions approximating 
$\varphi(x,\xi)=1$. Then, we write the Fourier transform of $u^\flat$ for $n\ne 0$:
\begin{align*}
\hat u^\flat(n,t)&\ds =\int_0^t \int_\R (\gamma |n|^{2\alpha}+\delta) e^{-s(-i a(\xi)\cdot n +\gamma |n|^{2\alpha}+\delta)} \widehat f(n,\xi,t-s)d\xi ds\\
&\ds = \int_0^t (\gamma |n|^{2\alpha}+\delta) e^{-s(\gamma |n|^{2\alpha}+\delta)} \mathcal G \widehat f(n,ns,t-s)ds,
\end{align*}
with, as in \eqref{defGG},
$$
 \mathcal G \widehat f(n,z,s)= \int_\R e^{-i a(\xi)\cdot z} \widehat f(n,\xi,s)d\xi.
$$
By Jensen's inequality, we have
$$
|\hat u^\flat(n,t)|^2 \le \int_0^t (\gamma |n|^{2\alpha}+\delta) e^{-s(\gamma |n|^{2\alpha}+\delta)} \left|\mathcal G \widehat f(n,ns,t-s)\right|^2 ds.
$$
Then, by simple manipulation and Lemma 2.4 in \cite{BouchutDesvillettes99}, we deduce the
fol\-lo\-wing sequence of inequality for $T\ge 0$ (recall that $\omega_n$ is defined by \refe{defomegan}):
\begin{align*}
\int_0^T |\hat u^\flat(n,t)|^2 dt &\ds \le \int_0^T \int_0^t (\gamma |n|^{2\alpha}+\delta) e^{-s(\gamma |n|^{2\alpha}+\delta)} \left|\mathcal G \widehat f(n,ns,t-s)\right|^2 ds dt\\
\\
&\ds \le \int_0^T \int_0^T (\gamma |n|^{2\alpha}+\delta) e^{-s(\gamma |n|^{2\alpha}+\delta)} \left|\mathcal G \widehat f(n,ns,t)\right|^2 ds dt\\
\\
&\ds \le \int_0^T \int_{\R^+} \omega_n e^{-\omega_n s} \left|\mathcal G \widehat f\left(n,\frac{n}{|n|}s,t\right)\right|^2 ds dt\\
\\
&\ds \le \int_0^T \int_{\R^+} \omega_n \frac1{1+\omega_n^2s^2} \left|\mathcal G \widehat f\left(n,\frac{n}{|n|}s,t\right)\right|^2 ds dt\\
\\
&\ds \le c\omega_n^b \int_0^T \left| \widehat f\left(n,t\right)\right|^2_{L^2_\xi} dt.
\end{align*}
We conclude by summing over $n\in\Z^N$, $n\ne 0$:
$$
\int_0^T\|u^\flat(t)\|_{H^{\left(\frac12-\alpha\right)b}}^2 dt \le c\gamma^{b}\int_0^T\left\| f(t)\right\|_{L^2_{x,\xi}}^2.
$$
Since $\left\| f(t)\right\|_{L^2_{x,\xi}}^2=\left\| u(t)\right\|_{L^1_{x}}$, we obtain
\begin{equation}\label{estim-uflat}
\int_0^T\|u^\flat(t)\|_{H^{\left(\frac12-\alpha\right)b}}^2 dt\le  c\gamma^{b}\int_0^T\left\| u(t)\right\|_{L^1_{x}}dt.
\end{equation}

\subsection{Estimate of $P$}\label{P}
Recall that $P$ is a random measure on $\T^N$ given by:
$$
\<P(t),\varphi\>=\sum_{k\geq 1} \int_{\T^N}\int_0^t \left(S_{A_{\gamma,\delta}}^*(t-s)\varphi\right)(x,u(s,x)) g_k(x) d\beta_k(s)dx,\, \varphi\in C(\T^N).
$$
Note that, for each $k\in\N$, 
\begin{equation*}
\begin{aligned}
\int_{\T^N}\int_0^t &\left(S_{A_{\gamma,\delta}}^*(t-s)\varphi\right)(x,u(s,x)) g_k(x) d\beta_k(s)dx\\
& = \int_{\T^N}\int_0^t e^{-\delta(t-s)}\left(e^{-B_\gamma(t-s)}\varphi\right)\big(x+a(u(s,x))(t-s)
\big) g_k(x) d\beta_k(s)dx\\
& = \int_{\T^N}\int_0^t e^{-\delta(t-s)}\varphi(x) e^{-B_\gamma(t-s)} 
\big(g_k(\cdot-a(u(\cdot,s))(t-s))\big) (x) d\beta_k(s)dx.
\end{aligned}
\end{equation*}
By \eqref{D0}, the mapping $x\mapsto g_k(x-a(u(x,s))(t-s))$ is a bounded function, therefore
by Lemma \ref{lemmahalfDelta} we can write for $\sigma\ge 0$:
$$
\|e^{B_\gamma(t-s)} \left(g_k(\cdot-a(u(\cdot,s))(t-s))\right)\|_{H^\sigma(\T^N)}\le 
c \left(\gamma (t-s) \right)^{-\frac\lambda{2\alpha}} \|g_k\|_{L^2(\T^N)}.
$$
We deduce that, for $\sigma\in [0,\alpha)$,  this defines
a function in $H^{\sigma}(\T^N)$ and
\begin{equation}\label{e57}
\begin{aligned}
\E&\left(\left\| \int_0^t e^{-\delta(t-s)}e^{-B_\gamma(t-s)} \left(g_k(\cdot-a(u(\cdot,s))(t-s))\right)d\beta_k(s) \right\|^2_{H^{\sigma}(\T^N)}\right)\\
&\le c \int_0^t e^{-2\delta(t-s)}\left(\gamma (t-s) \right)^{-\frac\sigma{\alpha}} ds\|g_k\|_{L^2(\T^N)}^2\\
&\le c \gamma^{-\frac\sigma\alpha}\delta^{\frac\sigma\alpha-1}\int_{\R^+}e^{2s}s^{-\frac\sigma\alpha}ds \|g_k\|_{L^2(\T^N)}^2
\end{aligned}
\end{equation}
Since the sum over $k$ of the right hand side of \eqref{e57} is finite, we may then conclude:
\begin{align*}
\E&\left(\left\| \sum_{k\ge 1}\int_0^t e^{B_\gamma(t-s)} \left(g_k(\cdot-a(u(\cdot,s))(t-s))\right)d\beta_k(s) \right\|^2_{H^{\sigma}(\T^N)}\right)\\
&\le c \gamma^{-\frac\sigma\alpha}\delta^{\frac\sigma\alpha-1}\left(\int_{\R^+}e^{-2s}s^{-\frac\sigma\alpha}ds\right) D_0 .
\end{align*}
This shows that in fact $P(t)$ is more regular than a measure. It is in the space $L^2(\Omega;H^\sigma(\T^N))$ for $\sigma < \alpha$ and 
\begin{equation}\label{estim-P}
\E\left(\left\| P(t)\right\|_{H^\sigma(\T^N)}^2 \right)\le 
c \gamma^{-\frac\sigma\alpha}\delta^{\frac\sigma\alpha-1}\left(\int_{\R^+}e^{-2s}s^{-\frac\sigma\alpha}ds\right) D_0. 
\end{equation}

\subsection{Bound on the kinetic measure}\label{sec:boundmm}

Before we estimate $Q$, we give an estimate on the kinetic measure.

\begin{lemma} \label{l8}Let $u\colon\T^N\times [0,T]\times\Omega\to\R$ be the solution to~\refe{stoSCL} with initial datum $u_0$. Then the measure $q:=m-\frac12\GG^2\delta_{u=\xi}$ satisfies
\begin{equation}\label{estimmsto}
\E\int_{\T^N\times[0,T]\times\R}\theta(\xi)d|q|(x,t,\xi)\leq D_0 \E\|\theta(u)\|_{L^1(\T^N\times[0,T])}
+ \E\int_{\T^N} \Theta(u_0(x))dx
\end{equation}
for all non-negative $\theta\in C_c(\R)$, where $\Theta(s)=\int_0^s\int_0^\sigma \theta(r)dr d\sigma$.
\label{lemmakineticmeasuresto}\end{lemma}

{\bf Proof:} Let $\theta\in C_c(\R)$ be non-negative. Set 
$$
\Theta(s)=\int_0^s\int_0^\sigma \theta(r)dr d\sigma.
$$
We test the kinetic formulation \refe{eq:kineticupreL1} with $\Theta'(\xi)$ and sum over $x\in\T^N$, $\xi\in\R$:
$$
\frac{d\;}{dt}\int_{\T^N}\E\Theta(u) dx=\frac12\E\int_{\T^N} \GG^2\theta(u)dx-\E\int_{\T^N\times\R}\theta(\xi) dm(x,t,\xi).
$$
Since $\Theta(0)=0$ and $\Theta(u)\geq 0$, we deduce  that
\begin{align*}
\E|\<|q|,\theta\>|&\leq \frac12\E\int_{\T^N\times[0,T]} \GG^2\theta(u)dx+\E\int_{\T^N\times[0,T]\times\R}\theta(\xi) dm(x,t,\xi),\\
&\leq \E\int_{\T^N\times[0,T]} \GG^2\theta(u)dx+ \E\int_{\T^N} \Theta(u_0(x))dx\\
& \leq D_0 \E\|\theta(u)\|_{L^1(\T^N\times[0,T])}+ \E\int_{\T^N} \Theta(u_0(x))dx. \hspace*{1cm}\qed
\end{align*}

\subsection{Estimate of $Q$}\label{Q}
 
 Recall that
 $$
 \<Q(t),\varphi\>=\int_{\T^N\times[0,t]\times\R} (t-s)a'(\xi)\cdot\nabla S_{A_{\gamma,\delta}}^*(t-s)\varphi dq(x,s,\xi).
 $$
 For $\eps>0$, we define the measure:
 $$
 \<Q_\eps(t),\varphi\>=\int_{\T^N\times[0,(t-\eps)\wedge 0]\times\R} (t-s)a'(\xi)\cdot\nabla S_{A_{\gamma,\delta}}^*(t-s)\varphi dq(x,s,\xi),
 $$ 
where $\varphi\in C(\T^N)$. We choose $\varphi$ depending also on $\omega\in \Omega$ and $t\in [0,T]$. We then have, for $\lambda\in(0,2]$,
\begin{align*}
\E&\int_0^T \langle (-\Delta)^{\frac\lambda2}Q_\eps(t),\varphi(t)\rangle dt\\
&= \E\int_\eps^T \int_{\T^N\times[0,t-\eps]\times\R} (t-s)a'(\xi)\cdot\nabla (-\Delta)^{\frac\lambda2} S_{A_{\gamma,\delta}}^*(t-s)\varphi(t) dq(x,s,\xi) dt.
\end{align*}
By  Lemma \ref{lemmahalfDelta}, for $p\ge 1$ and $p'$ its conjugate exponent, we have
\begin{align*}
\|(-\Delta)^{\frac\lambda2}\nabla S_{A_{\gamma,\delta}}^*(t-s)&\varphi(t)\|_{L^{\infty}(\T^N\times \R)}\\
&= e^{-\delta(t-s)}\|(-\Delta)^{\frac\lambda2}\nabla e^{B_\gamma(t-s)}\varphi(t)\|_{L^{\infty}(\T^N)}\\
& \le c(\gamma(t-s))^{-2+\mu_{N,\alpha,\lambda,p}} e^{-\delta(t-s)}\|\varphi(t)\|_{L^{p'}(\T^N)},
\end{align*} 
with 
\begin{equation}\label{defmuparam}
\mu_{N,\alpha,\lambda,p} =2- \frac{N+\lambda+1}{2\alpha}+\frac{N}{2p\alpha}.
\end{equation}
Assume  $\mu_{N,\alpha,\lambda,p}>0$. By the Fubini Theorem, we have then
\begin{align*}
 \E\int_0^T& \langle (-\Delta)^{\frac\lambda2} Q_\eps(t),\varphi(t)\rangle dt\\
& \le \frac{\kappa}{\gamma} \|\varphi\|_{L^\infty(\Omega
 \times (0,T);L^{p'}(\T^N))} \E\int_{\T^N\times[0,T]\times\R}
 |a'(\xi)| d|q|(x,s,\xi),
 \end{align*}
where 
 $$
\kappa=\sup_{s\in [0,T],\; T>0}\int_{s+\eps}^T(\gamma(t-s))^{-1+\mu}
 e^{-\delta(t-s)}dt \le \gamma^{-1+\mu}\delta^{-\mu}\int_{\R^+} \sigma^{\mu-1 }e^{-\sigma}d\sigma,
 $$
with $\mu=\mu_{N,\alpha,\lambda,p}$. It follows that for all $\eps>0$, $Q_\eps\in L^1(\Omega\times (0,T);W^{\lambda,p}(\T^N))$  with
\begin{align*}
\E\|Q_\eps\|&_{L^1((0,T);W^{\lambda,p}(\T^N))}\\
&\le C_{N,\alpha,\lambda,p}\gamma^{-2} 
 \left(\frac\gamma\delta\right)^{\mu_{N,\alpha,\lambda,p}}
\E\int_{\T^N\times[0,T]\times\R}
 |a'(\xi)| d|q|(x,s,\xi),
\end{align*}
for some constant $C_{N,\alpha,\lambda,p}$. We let $\eps\to 0$ and use Lemma \ref{lemmakineticmeasuresto} to obtain eventually the following estimate on $Q$: 
 \begin{equation}\label{estim-Q}
 \begin{aligned}
 &\E\|Q\|_{L^1((0,T);W^{\lambda,p}(\T^N))}\\
 &\le C_{N,\alpha,\lambda,p}\gamma^{-2} \left(\frac\gamma\delta\right)^{\mu_{N,\alpha,\lambda,p}}
\left(D_0 \E\| a'(u)\|_{L^1(\T^N\times (0,T))} + \E\int_{\T^N} \Theta(u_0) dx\right),
\end{aligned}
 \end{equation}
where $\Theta(s)=\int_0^s\int_0^\sigma |a'(r)|drd\sigma$, provided $\mu_{N,\alpha,\lambda,p}>0$.
 
\subsection{Conclusion and proof of the existence part in Theorem~\ref{th:maintheorem}}\label{conclusion}

Under the growth hypothesis~\eqref{aquad} (sub-linearity of $a'$), the bound \eqref{estim-Q} gives 
 \begin{equation*}\label{estim-QQ}
 \begin{aligned}
 &\E\|Q\|_{L^1((0,T);W^{\lambda,p}(\T^N))}\\
 &\le C_{N,\alpha,\lambda,p,D_0}\gamma^{-2} \left(\frac\gamma\delta\right)^{\mu_{N,\alpha,\lambda,p}}
\left(1+\E\|u\|_{L^1(\T^N\times (0,T))} + \E\|u_0\|_{L^3(\T^N)}^3\right).
\end{aligned}
 \end{equation*}
We choose $\gamma,\, \delta>0$ such that:
$$
C_{N,\alpha,\lambda,p,D_0}\gamma^{-2} \left(\frac\gamma\delta\right)^{\mu_{N,\alpha,\lambda,p}}\le \frac14 
$$ 
and obtain 
\begin{align} 
\E\|Q\|&_{L^1(0,T;W^{\lambda,p}(\T^N))}\nonumber\\
&\le \frac14\E\|u\|_{L^1(0,T;L^p(\T^1))} + C(N,\alpha,\lambda,p,\gamma,\delta,D_0)(\E\|u_0\|_{L^3(\T^N)}^3+1).\label{Qfinal}
\end{align}
Consider now the estimates \eqref{estim-u0}, \eqref{estim-uflat}, \eqref{estim-P} and assume that $H^{\left(\frac12-\alpha\right)b}$
and $H^\sigma$ are both imbedded in a given Sobolev space $W^{s,q}(\T^N)$ with $s> 0$ (the exponent $s$, $q$ will be determined later). Then \eqref{estim-u0}, \eqref{estim-uflat}, \eqref{estim-P} give
\begin{align*}
\E\big(\|u^0+u^\flat+P&\|_{L^2(0,T;W^{s,q}(\T^N))}^2\big)\nonumber\\
&\le C(\delta,\gamma,D_0)\left(1+\E\|u\|_{L^1(\T^N\times(0,T))}+\E\|u_0\|_{L^1(\T^1)}+T\right).
\end{align*}
By the Cauchy-Schwarz inequality and the Young inequality, we obtain
\begin{align}
\E\big(\|u^0+&u^\flat+P\|_{L^1(0,T;W^{s,q}(\T^N))}\big)\nonumber\\
&\le C(\delta,\gamma,D_0)\left(1+\E\|u_0\|_{L^1(\T^N)}+T\right)+\frac14\E\|u\|_{L^1(0,T;L^1(\T^N))}.\label{P0bfinal}
\end{align}
Assuming also that $W^{\lambda,p}(\T^N)$ is embedded in $W^{s,q}(\T^N)$, we deduce from \eqref{decompusto}, \eqref{Qfinal}, \eqref{P0bfinal} the estimate
\begin{align*}
\E\|u&\|_{L^1(0,T;W^{s,q}(\T^N))} \\
&\le 
C(N,\alpha,\lambda,p,\gamma,\delta,D_0)( \E\|u_0\|_{L^3(\T^1)}^3+1+T)+\frac12\E\|u\|_{L^1(0,T;L^1(\T^N))},
\end{align*}
and thus
\begin{equation}\label{FINALu}
\E\|u\|_{L^1(0,T;W^{s,q}(\T^N))} \le C(N,\alpha,\lambda,p,\gamma,\delta,D_0)( \E\|u_0\|_{L^3(\T^N)}^3+1+T).
\end{equation}
Recall that the estimates above are actually obtained for $u^\eta$, the solution of \eqref{parabolic} which converges to the solution of  \eqref{stoSCL} in $L^\infty(0,T;L^1(\Omega\times \T^1))$. Thanks 
to the closedness of balls of $L^1(\Omega\times (0,T);W^{\lambda,p}(\T^1))$ for the 
$L^\infty(0,T;L^1(\Omega\times \T^1))$ topology, these estimates actually hold for the true solution $u$. Then, by Krylov-Bogoliubov theorem and the compactness of the embedding of $W^{s,q}(\T^N)$ in $L^1(\T^N)$, existence of an invariance measure follows. Moreover, by the Sobolev embedding \eqref{SobolevEmbedding} below, this invariant measure is supported by $L^r(\T^N)$ for any $r$ such that $\frac{s}{N}-\frac1q\geq -\frac1r$.
\medskip

There remains to determine for which $s,q$ we do have \eqref{FINALu}. Let us recall that we have the injection
\begin{equation}\label{SobolevEmbedding}
W^{r,p}(\T^N)\hookrightarrow W^{s,q}(\T^N)
\end{equation}
if $r\geq s$ and $\frac{r}{N}-\frac1p\geq\frac{s}{N}-\frac1q$, where $\frac{s}{N}-\frac1q$ can be considered as the index of regularity of the Sobolev space 
$W^{s,q}(\T^N)$. Then \eqref{estim-u0}, \eqref{estim-uflat} and \eqref{estim-P} provide us with the indices
$$
\mathrm{ind}_0=\frac{\alpha+\left(\frac12-\alpha\right)b}{N}-\frac12,\;
\mathrm{ind}_\flat=\frac{\left(\frac12-\alpha\right)b}{N}-\frac12,\;
\mathrm{ind}_P=\frac{\alpha-}{N}-\frac12,
$$
where $\alpha-$ denotes any positive number $\sigma<\alpha$. The condition $\mu_{N,\alpha,\lambda,p}>0$ on the parameter $\mu_{N,\alpha,\lambda,p}$ (defined by \refe{defmuparam}) reads 
$$
\mathrm{ind}_Q:=\frac{\lambda}{N}-\frac1p<\mathrm{ind}_Q^+:=\frac{4\alpha}{N}-\frac{N+1}{N},
$$
where $\mathrm{ind}_Q$ is associated to \eqref{estim-Q}. Clearly, we have $\mathrm{ind}_0>\mathrm{ind}_\flat$. We look for a sobolev space $W^{s,q}(\T^N)$ of index lower than the minimum of $\mathrm{ind}_\flat$, $\mathrm{ind}_P$, $\mathrm{ind}_Q^+$, \textit{i.e.}
\begin{equation}\label{eq:indices}
s<\frac{N}{q}+\min\left(\left(\frac12-\alpha\right)b-\frac{N}{2},\alpha-\frac{N}{2},4\alpha-(N+1)\right),
\end{equation}
and satisfying also
\begin{align}
0<&s,\label{eq:ssobolev1}\\
s<&\min\left(\left(\frac12-\alpha\right)b,\alpha\right)\label{eq:ssobolev2}.
\end{align}
The constraint \eqref{eq:ssobolev2} is contained in \eqref{eq:indices} as soon as $q\geq 2$. The constraint \eqref{eq:ssobolev1}, together with \eqref{eq:indices} gives (for the maximal admissible value $\alpha=\frac12$) the condition $\frac1N+\frac1q>1$, which is compatible with $q\geq 2$ only if $N=1$. Therefore we treat separately the cases $N=1$ and $N>1$. 
\bigskip

{\bf Case  $N=1$:} We choose $q=2$ and then optimize the right-hand side of \refe{eq:indices}. The maximum is $\frac{b}{2(b+4)}$, reached for $\alpha=\frac{b+3}{2(b+4)}$. Then we choose $p=q$ and $s=\lambda\in\left(0,\frac{b}{2(b+4)}\right)$. We obtain that the invariant measure is supported by $L^r(\T^1)$ for any $r<2+\frac{b}{2}$.
\medskip

{\bf Case  $N>1$:} take $q$ satisfying $\frac1N+\frac1q=1+\frac{\eta}{N}$, $\eta>0$ small. We then have to consider the maximum of
$$
\min\left(\left(\frac12-\alpha\right)b,\alpha,4\alpha-2+\eta\right),
$$
for $\alpha\in[0,\frac12]$, which is $\frac{\eta b}{b+4}$, obtained for $\alpha=\frac12-\frac{\eta}{b+4}$. We choose then $p=q$, $\lambda=s\in\left(0,\frac{\eta b}{b+4}\right)$. We obtain that the invariant measure is supported by $L^r(\T^N)$ for any $r<\frac{N}{N-1}$. \qed

\section{Uniqueness of the invariant measure, er\-go\-di\-ci\-ty.}\label{sec:uniqueness}

We have proved in \eqref{FINALu} that there exists $q>1$ and $s>0$ such that:
\begin{equation}\label{e33}
\E\|u\|_{L^1(0,T;W^{s,q}(\T^N))} \le 
 \kappa_0( \E\|u_0\|_{L^{3}(\T^N)}^{3}+1+T)
\end{equation}
where $\kappa_0$ depends on $\lambda,p,D_0$.  Below, in Section~\ref{s4.1}, we prove that this 
estimate implies that any solution enters a ball of some fixed radius in finite time. Then in the case of 
subquadratic flux $A$, we 
show that if the driving noise is small for some time, then the solution becomes very small, \textit{cf.} Section~\ref{sec:noisesmall}. This 
smallness depends on the size of the initial data. These two facts are then shown to imply the second part of Theorem~\ref{th:maintheorem} (Section~\ref{sec:conclusion}).

\subsection{Time to enter a ball in $L^1(\T^N)$}\label{s4.1}

For further purposes, we need to consider two solutions $u^1$ and $u^2$ starting from 
$u^1_0$, $u^2_0$ deterministic and in $L^{3}(\T^N)$.
It is easy to generalize \eqref{e33} to two solutions on the interval $[t,t+T]$ for $t,\, T\ge 0$. This implies:
\begin{align}
\E\big(\int_t^{t+T}\|u^1(s)\|_{L^{1}(\T^N)}&+\|u^2(s)\|_{L^{1}(\T^N)} ds\big| \mathcal F_t\big)\nonumber\\
&\le \kappa_1( \|u^1(t)\|_{L^{3}(\T^N)}^{3}+\|u^2(t)\|_{L^{3}(\T^N)}^{3}+1+T)\label{e33bis}
 \end{align}


Moreover, we easily prove for $i=1,\, 2$, thanks to the It\^o formula:
\begin{align}
&\|u^i(t)\|^3_{L^3(\T^N)}\nonumber\\
&\le \|u^i_0\|^3_{L^3(\T^N)}  + 3\int_0^t ((u^i)^2(s),\Phi dW(s))_{L^2(\T^N)}
+3D_0\int_0^t \|u^i(s)\|_{L^1(\T^N)}ds.\label{e34}
\end{align}
It follows that
\begin{align*}
\E&\bigg(\int_t^{t+T}\|u^1(s)\|_{L^{1}(\T^N)}+\|u^2(s)\|_{L^{1}(\T^N)} ds\big| \mathcal F_t\bigg)\\
\le &\kappa_1\bigg[ \|u^1_0\|^3_{L^3(\T^N)}+\|u^2_0\|^3_{L^3(\T^N)} +3 D_0 \int_0^t \|u^1(s)\|_{L^1(\T^N)}+\|u^2(s)\|_{L^1(\T^N)}ds\\
&+ 3\int_0^t ((u^1(s))^2+(u^2(s))^2,\Phi dW(s))_{L^2(\T^N)}+1+T\bigg]
 \end{align*}

We now define recursively the sequences of deterministic  times $(t_k)_{k\ge 0}$ and 
$(r_k)_{k\ge 0}$ by
$$
\begin{array}{l}
\ds t_0=0,\\
\\
\ds t_{k+1}=t_k+r_k,
\end{array}
$$
where $(r_k)_{k\ge 0}$ will be chosen below. We also define the events:
$$
A_k=\left\{ \inf_{s\in [t_\ell,t_{\ell+1}]}\big(\|u^1(s)\|_{L^{1}(\T^N)}+\|u^2(s)\|_{L^{1}(\T^N)}\big)\ge 2\kappa_1,\;
\ell =0,\dots, k-1\right\}. 
$$
Then, for all $k\ge 0$, 
\begin{align}
&\P\left(\inf_{s\in [t_k,t_{k+1}]}\big(\|u^1(s)\|_{L^{1}(\T^N)}+\|u^2(s)\|_{L^{1}(\T^N)}\big)\ge 2\kappa_1 \bigg| \mathcal F_{t_k}\right)\nonumber\\
\le &\P\left(\frac1{r_k} \int_{t_k}^{t_k+r_k}\|u^1(s)\|_{L^{1}(\T^N)}+\|u^2(s)\|_{L^{1}(\T^N)} ds\ge 2\kappa_1 \bigg| \mathcal F_{t_k}\right)\nonumber\\
\le &\frac1{2r_k}\bigg[ \|u^1_0\|^3_{L^3(\T^N)}+\|u^2_0\|^3_{L^3(\T^N)} \nonumber\\
&+3 D_0 \int_0^{t_k} \|u^1(s)\|_{L^1(\T^N)}+\|u^2(s)\|_{L^1(\T^N)}ds +1\bigg]\nonumber\\
 &+\frac12 + \frac3{2r_k}\int_0^{t_k} ((u^1(s))^2+(u^2(s))^2,\Phi dW(s))_{L^2(\T^N)}.\label{ineqAk}
\end{align}
We will choose $r_k$ so that the following inequality is satisfied for all $k\ge 0$:
$$
\frac1{2r_k}( \|u^1_0\|^3_{L^3(\T^N)}+\|u^2_0\|^3_{L^3(\T^N)}  +1)\le \frac18.
$$
Let us then multiply \eqref{ineqAk} by $\carac_{A_k}$ and take the expectation to obtain
$$
\begin{array}{ll}
\ds
\P\left(A_{k+1}\right) \le &\ds \frac58\P\left(A_{k}\right)+ \frac{3D_0}{2r_k}\E\left(\int_0^{t_k} \|u^1(s)\|_{L^1(\T^N)} +\|u^2(s)\|_{L^1(\T^N)}ds\carac_{A_k}\right)\\
&\ds
+
\frac3{2r_k}\E\left(\int_0^{t_k} ((u^1(s))^2+(u^2(s))^2,\Phi dW(s))_{L^2(\T^N)}\carac_{A_k}\right).
\end{array}
$$
By the Cauchy-Schwarz inequality, It\^o isometry and the bound \eqref{D0} on the covariance of the noise, we have:
\begin{align*}
\E\big(\int_0^t& ((u^1(s))^2+(u^2(s))^2,\Phi dW(s))_{L^2(\T^N)}\carac_{A_k}\big)\\
& \le (2D_0)^{1/2} \left( \E\int_0^t \|u^1(s)\|_{L^2(\T^N)}^4+\|u^2(s)\|_{L^2(\T^N)}^4ds\right)^{1/2}\P(A_k)^{1/2}.
\end{align*}
By It\^o formula applied to $\|u^i(t)\|_{L^2(\T^N}^2$ and 
$\|u^i(t)\|_{L^2(\T^N}^4$, $i=1,2$, we obtain the bound:
$$
\E\left(\|u^i(t)\|_{L^2(\T^N)}^2\right)\le \E\left(\|u^i_0\|_{L^2(\T^N)}^2\right)+D_0t
$$
and
$$
\E\left(\|u^i(t)\|_{L^2(\T^N)}^4\right)\le \E\left(\|u^i_0\|_{L^2(\T^N)}^4+6 D_0\int_0^t\|u^i(s)\|_{L^2(\T^N)}^2ds\right).
$$
Thus, integrating in time, we obtain a bound
$$
\E\int_0^{t_k} \|u^1(s)\|^4_{L^2(\T^N)}+\|u^2(s)\|^4_{L^2(\T^N)}ds
\le 2[\|u^1_0\|^4_{L^2(\T^N)}+\|u^2_0\|^4_{L^2(\T^N)}]t_k + 24 D_0^2t_k^3.
$$
Similarly, we have the estimate
$$
\E\left|\int_0^{t_k} \|u^1(s)\|_{L^1(\T^N)}+\|u^2(s)\|_{L^1(\T^N)}ds\right|^2
\le t_k^2[\|u^1_0\|^2_{L^2(\T^N)}+\|u^2_0\|^2_{L^2(\T^N)}]+D_0 t_k^3.
$$
It follows that
$$
\P\left(A_{k+1}\right) \le \frac34\P\left(A_{k}\right)+ 
\frac{C}{r_k^2}\left([\|u^1_0\|^4_{L^2(\T^N)}+\|u^2_0\|^4_{L^2(\T^N)}]t_k +  t_k^4 +1 \right),
$$
where the constant $C$ can be written explicitely in term of $D_0$.
We choose $r_k$ so that:
$$
\frac{C}{r_k^2}\left([\|u^1_0\|^4_{L^2(\T^N)}+\|u^2_0\|^4_{L^2(\T^N)}]t_k +  t_k^4 +1 \right)\le \left(\frac34\right)^k,
$$
and obtain:
$
\P\left(A_{k+1}\right) \le \frac34\P\left(A_{k}\right)+ \left(\frac34\right)^k,
$
which gives therefore 
$$
\P\left(A_{k}\right)\le k\left(\frac34\right)^{k-1}.
$$
By Borel-Cantelli Lemma, we deduce that
$$
k_0=\inf \{  k\ge 0\; |\; \inf_{s\in [t_\ell,t_{\ell+1}]}\|u^1(s)\|_{L^{1}(\T^N)}+\|u^2(s)\|_{L^{1}(\T^N)}\le 2\kappa_1  \}
$$
is almost surely finite. We then define the stopping time 
$$
\tau^{u^1_0,u^2_0}=\inf\{t\ge 0\;|\;  \|u^1(t)\|_{L^{1}(\T^N)}+\|u^2(t)\|_{L^{1}(\T^N)} \le 2\kappa_1 \}.
$$
Clearly $\tau^{u^1_0,u^2_0}\le t_{k_0+1}$ so that $\tau^{u^1_0,u^2_0}<\infty$ almost surely.
It follows that for $T>0$ the following stopping times are also almost surely finite:
$$
\tau_\ell=\inf\{t\ge \tau_{\ell-1}+T\;|\; \|u^1(t)\|_{L^{1}(\T^N)}+\|u^2(t)\|_{L^{1}(\T^N)} \le 2\kappa_1 \},\; \tau_0=0.
$$

\subsection{The solution is small if the noise is small}\label{sec:noisesmall}


\begin{proposition}\label{p9} Assume that $a$ satisfies \eqref{linear}. Then, for any $\eps>0$, there exists $T>0$ and $\eta >0$ such that:
$$
\frac1T\int_0^T \|u(s)\|_{L^1(\T^N)}ds \le \frac\eps2
$$
if 
$$
 \|u(0)\|_{L^1(\T^N)}\le 2\kappa_1
\mbox{ and } \sup_{t\in[0,T]}\|W\|_{W^{1,\infty}(\T^N)	}\le \eta.
$$
\end{proposition}

{\bf Proof:} We first take $\tilde u_0\in L^2(\T^N)$ such that 
$$
\|u_0-\tilde u_0\|_{L^1(\T^N)}\le \frac{\eps}{8}, \; \|\tilde u_0\|_{L^2(\T^N)}\le C \kappa_1 \eps^{-N/2} .
$$
It is easy to see that this can be achieved by taking 
$\tilde u_0=\rho_\eps * u_0$ for some regularizing kernel $(\rho_\eps)_{\eps>0}$.
By \eqref{contractionL1}, we know that
$$
\|u(t)-\tilde u(t)\|_{L^1(\T^N)}\le \frac{\eps}{8}, \; t\ge 0,
$$
where $\tilde u$ is the solution starting from $\tilde u_0$. We set $v=\tilde u-W$ so that $v$ is a kinetic solution to 
$$
\partial_t v+\div A (v+W)=0
$$
and $g=\carac_{v(x,t)>\xi}-\carac_{0>\xi}$ satisfies:
\begin{equation} \label{e35}
\partial_t g+a(\xi)\cdot \nabla g= \partial_\xi \mathrm{n} +(a(\xi)-a(\xi+W))\cdot \nabla g -a(\xi+W)\cdot \nabla W 
\delta_{v=\xi}.
\end{equation}
The kinetic measure $\mathrm{n}$ is easily related to the kinetic measure in the equation satisfied by 
$f=\carac_{u>\xi}-\carac_{0>\xi}$. Also since both $u$ and $W$ have a zero spatial average, so does $v$.
\medskip

We use again $B_\gamma$ and $A_{\gamma,\delta}$ introduced in section  \ref{sec:decu}. We 
now take $\alpha=\frac12$. Then \eqref{e35} rewrites:
\begin{equation} \label{e36}
\partial_t g+A_{\gamma,\delta} g= (B_\gamma+\delta \mbox{Id})g+ \partial_\xi \mathrm{n} +(a(\xi)-a(\xi+W))\cdot \nabla g -a(\xi+W)\cdot \nabla W 
\delta_{v=\xi}.
\end{equation}
By solving \eqref{e36} and summing over $\xi$, we are led to the following decomposition of $v$:
$$
v=v^0+v^\flat+v^\#+P^W+N^W,
$$
where
\begin{align*}
v^0(t)&:=\int_\R S_{A_{\gamma,\delta}}(t)g(0,\xi)d\xi,\\
v^\flat(t)&:=\int_\R \int_0^t S_{A_{\gamma,\delta}}(s)\, (B_\gamma +\delta \mbox{Id})g(t-s,\xi)ds d\xi,\\
v^\#(t)&:=\int_\R \int_0^t S_{A_{\gamma,\delta}}(s)(a(\xi)-a(\xi+W))\cdot \nabla g (t-s,\xi)\, ds d\xi,\\
\<P^W(t),\varphi\>&:=-\int_{\T^N\times[0,t]}\hspace{-.5cm}  \left(a(v+W)\cdot \nabla W \right)(x,s)\, 
\left(S_{A_\gamma,\delta}^*(t-s)\varphi\right)(x,v(x,s)) \,ds dx,\\
\<N^W(t),\varphi\>&:=\int_{\T^N\times[0,t]\times\R} \hspace{-.5cm}(t-s)a'(\xi)\cdot\nabla S_{A_{\gamma},\delta}^*(t-s)\varphi d\mathrm{n}(x,s,\xi),
\end{align*}
where 
$\<\cdot,\cdot\>$ is the duality product between the space of finite Borel measures on $\T^N$ and $C(\T^N)$.
\medskip

Reproducing the argument of Section \ref{u0} and \ref{uflat}, we obtain the estimates:
\begin{equation}\label{estimv0}
\int_0^T\|v^0(t)\|_{H^{\frac{1}2}(\T^N)}^2 dt \le c\,\gamma^{b-1}\left\| u_0\right\|_{L^1(\T^N)}
\end{equation}
and
\begin{equation}\label{estimvflat}
\int_0^T\|v^\flat(t)\|_{L^2(\T^N)}^2 dt \le  c\,\gamma^{b}\int_0^T\left\| v(t)\right\|_{L^1(\T^N)}dt.
\end{equation}
We estimate $v^\#$ similarly. We take the Fourier transform in $x$ (denoted by $\mathcal{F}$ here) and integrate in time. Note that
\begin{align*}
{\mathcal F}&\left[(a(\cdot)-a(\cdot+W))\cdot \nabla g\right](n,\xi,t)\\
&=2\pi i {\mathcal F}\left[(a(\cdot)-a(\cdot+W)) g\right](n,\xi,t)\cdot n
+{\mathcal F}\left[(a'(\cdot+W)\cdot \nabla W)g \right](n,\xi,t)
\end{align*} 
and thus
\begin{align*}
&\int_0^T |\hat v^\#(n,t)|^2dt \\
\leq& c   \int_0^T \left| \int_\R \int_0^t e^{-i(a(\xi)\cdot n +\gamma  |n|+\delta)s}
{\mathcal F}\left[(a(\cdot)-a(\cdot+W)) g\right](n,\xi,t-s)\cdot n \, ds d\xi \right|^2 dt\\
+&c \int_0^T \left| \int_\R \int_0^t e^{-i(a(\xi)\cdot n +\gamma  |n|+\delta)s}
{\mathcal F}\left[(a'(\cdot+W)\cdot \nabla W)g \right](n,\xi,t-s) \, ds d\xi \right|^2 dt.
\end{align*}
Denote by $A_1$ and $A_2$ the two terms on the right hand side. $A_1$ is of the form
$$
A_1=\int_0^T\left| \int_0^t e^{-(\gamma |n|+\delta) s} n\cdot {\mathcal G}h(n,ns,t-s)ds\right|^2dt,
$$
where 
$$
h(n,\xi,t)=\F\left[\left(a(\cdot)-a(\cdot+W)\right)g\right](n,\xi,t)
$$ 
and the integral operator $\mathcal{G}$ has been defined in \eqref{defGG}. We proceed as in Section \ref{u0} and \ref{uflat} (in particular we use Lemma 2.4 of \cite{BouchutDesvillettes99} for the estimate of the oscillatory integral $\mathcal{G}h$) to obtain
$$
A_1\leq c \gamma^{-2+b} \int_0^T \left\| h(n,\xi,t)\right\|^2_{L^2(\xi)}dt,
$$
and, similarly,
$$
A_2\le \frac{c \gamma^{-2+b}}{|n|} \int_0^T \|k(n,\xi,t)\|_{L^2_\xi}^2 dt
$$
with $k(n,\xi,t)={\mathcal F}\left[(a'(\cdot+W)\cdot \nabla W)g \right](n,\xi,t)$.
Summing over $n\ne 0$, we get therefore 
\begin{align*}
\int_0^T &\| v^{\#,1}(t)\|_{L^2(\T^N)}^2dt\\
&\le C\gamma^{-2+b} \int_0^T \|\left(a(\cdot)-a(\cdot+W)\right)g\|_{L^2_{x,\xi}}^2
+\|(a'(\cdot+W)\cdot \nabla W)g\|_{L^2_{x,\xi}}^2dt,
\end{align*}
where $v^{\#,1}=v^\#-\int_{\T^N}v^\#dx$. We then use the hypothesis \eqref{linear}, \textit{i.e.} $a$ sublinear, and the ``smallness" of the noise (actually no smallness assumption on $\eta$ has been done up to now) to write
$$
|a(\xi)-a(\xi+W)|\le c\,|W|\le c\,\eta,\quad |a'(\xi+W)\cdot\nabla W |\le c\,\eta
$$
and deduce
\begin{equation}\label{estimvsharp}
\int_0^T \| v^{\#,1}(t)\|_{L^2(\T^N)}^2dt\le C \eta^2 \gamma^{-2+b}\int_0^T \|v\|_{L^{1}(\T^N)}dt
\end{equation}
since $\|g\|_{L^2_{x,\xi}}^2=\|v\|_{L^{1}(\T^N)}$. The total mass $\int_{\T^N}v^\#dx$ will be estimated later, once we have finished to estimate $P^W$ and $N^W$.
Regarding $P^W$, we have, by Hypothesis~\eqref{linear},
\begin{align*}
\<P^W(t),\varphi\>\le C_{} \int_0^t (1+&\|v(s)\|_{L^1(\T^N)}+\|W(s)\|_{L^1(\T^N)}) \\
&\times \|\nabla W(s)\|_{L^\infty(\T^N)} 
\|S_{A_\gamma,\delta}^*(t-s)\varphi\|_{L^\infty(\T^N)} \,ds.
\end{align*}
Therefore $P^W$ is more regular than a measure and Lemma \ref{lemmahalfDelta}
implies, for $\eta\le 1$,
$$
\|P^W(t)\|_{L^1(\T^N)}\le C \eta  \int_0^t(1+\|v(s)\|_{L^1(\T^N)})e^{-\delta (t-s)}ds.
$$
In particular, we have
\begin{equation}\label{estimPW}
\int_0^T \|P^W(t)\|_{L^1(\T^N)}dt 
\le C\frac{\eta}{\delta} \int_0^T (1+\|v(s)\|_{L^1(\T^N)}) ds.
\end{equation}
We finally estimate the last term containing the measure $\rm n$. We write, using 
again Lemma \ref{lemmahalfDelta} and the hypothesis~\eqref{linear},
\begin{align*}
\<N^W(t),\varphi\>&\ =\int_{\T^N\times[0,t]\times\R} (t-s)a'(\xi)\cdot\nabla S_{A_{\gamma}}^*(t-s)\varphi d\mathrm{n}(x,s,\xi)\\
&\ds \le c\,\|\varphi\|_{L^\infty(\T^N)} \int_{\T^N\times[0,t]\times\R} 
e^{-\delta (t-s)}  d|\mathrm{n}|(x,s,\xi).
\end{align*}
Again, we may prove that $N^W$ is more regular than a measure and deduce:
\begin{equation}\label{estimN0}
\int_0^T \|N^W(t)\|_{L^1(\T^N)} dt\le \frac{C}\delta \int_{\T^N\times[0,T]\times\R}   d|\mathrm{n}|(x,s,\xi).
\end{equation}
To complete our estimate, it remains to evaluate the mass of the measure $\mathrm{n}$. We proceed as in Lemma \ref{l8} and test equation \eqref{e36} against $\xi$ to obtain:
\begin{align*}
&\frac12\|v(t)\|^2_{L^2(\T^N)} + |\rm n| ([0,t]\times\T^N\times \R) \\
\le&\frac12  \|u_0\|^2_{L^2(\T^N)}+\left|\int_{[0,t]\times\T^N\times \R} \xi (a(\xi)-a(\xi+W))\cdot \nabla g \, dx \, d\xi \, ds\right|\\
+&\left|\int_{[0,t]\times\T^N} v(x,s) a(v(x,s) +W(x,s))\cdot \nabla W(x,s)dxds\right|.\\
\end{align*}
We use an integration by parts in the second term. By \eqref{linear} we obtain after easy manipulations:
$$
\frac12\|v(t)\|^2_{L^2(\T^N)} + |\mathrm{n}| ([0,t]\times\T^N\times \R) \le \frac12  \|u_0\|^2_{L^2(\T^N)}
+ \eta C\int_0^t (\|v(s)\|^2_{L^2(\T^N)}+1)ds.
$$
Gronwall lemma then gives
$
\|v(t)\|^2_{L^2(\T^1)}\le e^{\eta C t}(\|u_0\|^2_{L^2(\T^1)} +C),
$
and therefore
$$
|\mathrm{n}| ([0,T]\times\T^N\times \R) \le e^{\eta C T}(\|u_0\|^2_{L^2(\T^N)} +C).
$$
Since $\|\tilde u_0\|_{L^2(\T^N)}\leq C\kappa_1\eps^{-N/2}$, it follows by \eqref{estimN0} that
\begin{equation}\label{estimNW}
\int_0^T \|N^W(t)\|_{L^1(\T^N)} dt\le  \frac{C }\delta e^{\eta C T}(\kappa_1^2\eps^{-N} +1).
\end{equation}
Finally, since $v$, $v^0$ and $v^\flat$ have zero spatial averages, we have
$$
\int_{\T^N}v^\# +P^W+N^W dx=0. 
$$
Consequently we can estimate the total mass of $v^\#$ by 
$$
\left|\int_{\T^N}v^\# dx\right| \le  \|P^W\|_{L^1(\T^N)}+ \|N^W\|_{L^1(\T^N)}.
$$
We may now gather all the estimates obtained, \textit{i.e.} \eqref{estimv0}, \eqref{estimvflat}, \eqref{estimvsharp}, \eqref{estimPW}, \eqref{estimNW}, and deduce:
\begin{align*}
\int_0^T\|v^0(t)&+v^\flat(t)+v^{\#,1}(t)\|_{L^1(\T^1)}dt \\
&\le T^{1/2}\left(\int_0^T \|v^0(t)\|_{L^2(\T^1)}^2
+\|v^\flat(t)\|_{L^2(\T^1)}^2+\|v^{\#,1}(t)\|_{L^2(\T^1)}^2 dt \right)^{1/2}\\
& \le CT^{1/2}\left(\gamma^{b-1}\kappa_1 + (\gamma^b+\eta^2\gamma^{-2+b})\int_0^T\|v(t)\|_{L^1(\T^1)}dt\right)^{1/2}\\
& \le \frac14\int_0^T\|v(t)\|_{L^1(\T^1)}dt + CT (\gamma^b+\eta^2\gamma^{-2+b})+C(T\gamma^{b-1}\kappa_1)^{1/2} 
\end{align*}
and
\begin{align*}
\int_0^T \bigg\|\int_{\T^N}v^\#(x)dx+&P^W(t)+N(t)\bigg\|_{L^1(\T^N)}\\
&\le C\left( \frac\eta\delta  \int_0^T(1+\|v(t)\|_{L^1(\T^N)})dt+\frac{1 }\delta e^{\eta C T}(\kappa_1^2\eps^{-N} +1)\right).
\end{align*}
Let $r>0$ be a ratio that we will fix later. We choose (in that order), $\gamma$, $T$, $\delta$ such that
$$
C\gamma^b\leq r\eps,\quad C\left(\frac{\gamma^{b-1}\kappa_1}{T}\right)^{1/2} \leq r\eps,\quad \frac{C}\delta e^{C T}(\kappa_1^2\eps^{-N} +1)\leq r\eps,
$$
and then $\eta$ small enough so that
$$
C\frac\eta\delta\leq \min(r\eps,\frac14),\quad C\eta^2\gamma^{-2+b}\leq r\eps.
$$
With those choices, and for $r=\frac{1}{40}$, we obtain, if furthermore $\eta\le \frac{\eps}{8}$,
$$
\frac1T\int_0^T\|v(t)\|_{L^1(\T^N)}dt\le \frac\eps{4},\quad \frac1T\int_0^T\|\tilde u(t)\|_{L^1(\T^N)}dt\le \frac{3
\eps}{8}
$$
and
$$
\frac1T\int_0^T\|u(t)\|_{L^1(\T^N)}dt\le \frac\eps{2}<\eps.
$$
This is the desired conclusion. \qed

\subsection{Conclusion}\label{sec:conclusion}

We assume that \eqref{linear} holds.
Let $u_0^1,\; u_0^2$ be in $L^1(\T^N)$. Let $\eps>0$, we take $\tilde u_0^1,\; \tilde u_0^2$ in $L^3(\T^N)$ such that $\|u_0^i-\tilde u_0^i\|_{L^1(\T^N)}\le \frac\eps4$. . We denote by
$u^1,\; u\;\tilde u^1,\;\tilde u^2$ the corresponding solutions. We associate to   
$\tilde u_0^1,\; \tilde u_0^2$ the sequence of stopping times constructed in section 
\ref{s4.1}. We choose $T$ and $\eta$ given by  Proposition \ref{p9} and obtain thanks to the $L^1$-contraction~\eqref{contractionL1},
\begin{align*}
\P\bigg(\frac1T\int_{\tau_\ell}^{\tau_\ell+T}& \|u^1(s)-u^2(s)\|_{L^1(\T^N)}ds  \le \eps\; \bigg|\; \mathcal F_{\tau_\ell}\bigg)\\
&\ge \P\bigg(\frac1T\int_{\tau_\ell}^{\tau_\ell+T} \|\tilde u^1(s)-\tilde u^2(s)\|_{L^1(\T^N)}ds \ \le \frac\eps2\; \bigg|\; \mathcal F_{\tau_\ell}\bigg)\\
&\ge\P\bigg( \sup_{[\tau_\ell,\tau_\ell+T]}\|W(t)-W(\tau_\ell)\|_{W^{1,\infty}(\T^N)}\le \eta\; \bigg|\; \mathcal F_{\tau_\ell}\bigg).
\end{align*}
By the strong Markov property, the right hand side is non random and independent on $\ell$. It is 
clearly positive. We denote it by $\lambda$. For $\ell_0, k\in\N$ we then have
$$
\P\left(\frac1T\int_{\tau_\ell}^{\tau_\ell+T} \|u^1(s)-u^2(s)\|_{L^1(\T^N)}ds  \ge \eps,\; 
\mbox{ for } \ell=\ell_0,\dots,\ell_0+k\right) \le (1-\lambda)^k,
$$
and therefore 
\begin{align}
\P\bigg(\lim_{\ell \to \infty}\frac1T&\int_{\tau_\ell}^{\tau_\ell+T} \|u^1(s)-u^2(s)\|_{L^1(\T^N)}ds  \ge \eps\bigg)\nonumber\\
&=\P\bigg(\exists \ell_0\in\N\,;\frac1T\int_{\tau_\ell}^{\tau_\ell+T} \|u^1(s)-u^2(s)\|_{L^1(\T^N)}ds  \ge \eps,\; \mbox{ for } \ell\ge\ell_0\bigg)\nonumber\\
&=0.\label{trajmeet0}
\end{align}
%
%
Note the limit exists since, by \eqref{contractionL1}, $t\mapsto\|u^1(t)-u^2(t)\|_{L^1(\T^N}$ is a.s. non-increasing.
This latter property is again used to deduce from \eqref{trajmeet0} that 
$$
\P\left(\lim_{t \to \infty} \|u^1(t)-u^2(t)\|_{L^1(\T^N)}  \ge \eps\right)=0.
$$
We have thus proved:
$$
\lim_{t \to \infty} \|u^1(t)-u^2(t)\|_{L^1(\T^N)}=0,\; a.s.
$$
This easily implies the second statement (uniqueness part) of Theorem~\ref{th:maintheorem}.

\begin{appendix}
\section{Solutions in $L^1$}\label{app:L1}

In this section, we develop the $L^1$-theory for the Cauchy Problem associated to \refe{stoSCL}. In the deterministic framework, this has been done in B\'enilan, Carrillo, Wittbold \cite{BenilanCarrilloWittbold00} and in \cite{ChenPerthame03}. In \cite{BenilanCarrilloWittbold00}, a concept of {\it renormalized entropy solution} is introduced. In \cite{ChenPerthame03} it is noticed that the kinetic formulation for scalar (parabolic degenerate here) conservation laws can be expanded from the $L^\infty$ to the $L^1$ framework with quite minor adaptations, and this is the approach we will follow here: although, for $u\in L^1$, it may be impossible to give a sense to $\div(A(u))$, still we will show that the problem is well-posed. Note that the basic reason for this is that (on the whole space at least) non-uniqueness for \refe{stoSCL}, in the case $g_k=0$, requires some growth at infinity: this fact is illustrated by the counter-examples constructed in Goritskii, Panov \cite{GoritskiiPanov02} in particular.

\subsection{Generalized solutions}

In \cite{DebusscheVovelle10}, \cite{DebusscheVovelle10revised}, the notion of solution is extended to a notion of generalized solution. The notion of solution given by Definition~\ref{defkineticsolL1} can be extended accordingly. The proof of Theorem \ref{th:ResolutionL1} is a rather straightforward extension of the $L^\infty$ result. Only the decay of the kinetic measure and Lemma~7 in \cite{DebusscheVovelle10revised}, needed to prove time continuity, have to be revised. We will do this below in sections \ref{decay} and  \ref{app:secuniqueness}
\medskip

\begin{definition}[Young measure] Let $(X,\lambda)$ be a finite measure space. Let $\mathcal{P}_1(\R)$ denote the set of probability measures on $\R$. We say that a map $\nu\colon X\to\mathcal{P}_1(\R)$ is a Young measure on $X$ if, for all $\phi\in C_b(\R)$, the map $z\mapsto \nu_z(\phi)$ from $X$ to $\R$ is measurable. We say that a Young measure $\nu$ vanishes at infinity if 
\begin{equation}
\int_X\int_\R |\xi| d\nu_z(\xi)d\lambda(z)<+\infty.
\label{nuvanish}\end{equation}
\end{definition}

\begin{definition}[Kinetic function] Let $(X,\lambda)$ be a finite measure space. A measurable function $f\colon X\times\R\to[0,1]$ is said to be a kinetic function if there exists a Young measure $\nu$ on $X$ that vanishes at infinity such that, for $\lambda$-a.e. $z\in X$, for all $\xi\in\R$,
\begin{equation*}
f(z,\xi)=\nu_{z}(\xi,+\infty).
\end{equation*}
We say that $f$ is an {\rm equilibrium} if there exists a measurable function $u\colon X\to\R$ such that $f(z,\xi)=\mathbf{1}_{z>\xi}$ a.e., or, equivalently, $\nu_z=\delta_{u(z)}$ for a.e. $z\in X$.
\end{definition}

\begin{definition}[Generalized solution] Let $f_0\colon\Omega\times\T^N\times\R\to[0,1]$ be a kinetic function. A measurable function $f\colon\Omega\times\T^N\times[0,T]\times\R\to[0,1]$ is said to be a generalized solution to~\refe{stoSCL} with initial datum $f_0$ if $(f(t))$ is predictable and is a kinetic function such that $\nu:=-\partial_\xi f$ satisfies 
\begin{equation}
\E\left(\esssup_{t\in[0,T]}\int_{\T^N}\int_\R|\xi| d\nu_{x,t}(\xi) dx\right)<+\infty,
\label{eq:integrabilityf}\end{equation}
and such that there exists a kinetic measure $m$ such that for all $\varphi\in C^1_c(\T^N\times[0,T)\times\R)$, 
\begin{equation}
\label{eq:kineticfpre}
\begin{array}{l}
\ds \int_0^T\<f(t),\partial_t \varphi(t)\>dt+\<f_0,\varphi(0)\>
+\int_0^T \<f(t),a(\xi)\cdot\nabla\varphi(t)\>dt\\
\ds =-\sum_{k\geq 1}\int_0^T\int_{\T^N}\int_\R g_k(x,\xi)\varphi(x,t,\xi)d\nu_{x,t}(\xi)dxd\beta_k(t)\\
\ds -\frac{1}{2}\int_0^T\int_{\T^N}\int_\R \partial_\xi\varphi(x,t,\xi)\GG^2(x,\xi)d\nu_{(x,t)}(\xi) dx dt+m(\partial_\xi\varphi), \mbox{ a.s.}
\end{array}
\end{equation}
\label{d4L1}\end{definition}

What differ mostly between  Definition~\ref{def:kineticmeasureL1} of kinetic measure and the corresponding definition in \cite{DebusscheVovelle10} is the condition at infinity {\it 2}. This condition is for example crucial to show the equivalence between Definition~\ref{defkineticsolL1} and the notion of renormalized solution (what we will not do here since we will not need the latter). As we will see, this decay condition is related to the decay of $k\mapsto (u-k)^\pm$. We explore this link and then use it to prove that solutions in the sense of Definition~\ref{defkineticsolL1} lead to a well-posed problem.

\subsection{Decay of the kinetic measure}\label{decay}

\begin{proposition} Let $u_{0}\in L^1(\T^N)$. Let a measurable function $u\colon\T^N\times [0,T]\times\Omega\to\R$ be a solution to~\refe{stoSCL} according to Definition~\ref{defkineticsolL1}. Let $T>0$. There exists a decreasing function $\eps\colon\R_+\to\R_+$ with $\lim_{k\to+\infty}\eps(k)=0$ depending on $T$ and on the functions
$$
k\mapsto \|(u_0-k)^+\|_{L^1(\T^N)},\quad k\mapsto \|(u_0-k)^-\|_{L^1(\T^N)}
$$
only such that, for all $k\geq 1$, 
$$
\E\left(\esssup_{t\in[0,T]}\|(u(t)-k)^\pm\|_{L^1(\T^N)}\right)+\E m(A_k)\leq\eps(k),
$$
where $A_k=\T^N\times[0,T]\times\{\xi\in\R,k\leq|\xi|\leq k+1\}$.
\label{prop:decaym}\end{proposition}

{\bf Proof.} {\bf Step 1.} For $k\ge 0$, set
$$
\theta_k(u)=\mathbf{1}_{k<u<k+1},\quad \Theta_k(u)=\int_0^{u}\int_0^r\theta_k(s) ds dr.
$$
Let $\gamma\in C^1_c([0,T))$ be non-negative and satisfy $\gamma(0)=1$, $\gamma'\leq 0$. After a preliminary step of approximation that uses the monotone convergence Theorem, we take $\varphi(x,t,\xi)=\gamma(t)\Theta_k'(\xi)$ in \refe{eq:kineticupreL1} to obtain
\begin{multline}
\int_0^T\int_{\T^N}\Theta_k(u(x,t))|\gamma'(t)| dx dt+\int_{A_k^+} \gamma(t) dm(x,t,\xi)\\
=\frac{1}{2}\E\int_0^T\int_{\T^N} \gamma(t)\GG^2(x)\theta_k(u(x,t)) dx dt+\int_{\T^N}\Theta_k(u_0(x))dx\\
+\sum_{j\geq 1}\int_0^T\int_{\T^N}g_j(x)\Theta_k'(u(x,t))\gamma(t) dx d\beta_k(t),
\label{eq:kineticupreL1gamma}\end{multline}
where $A_k^+=\T^N\times[0,T]\times\{\xi\in\R,k\leq \xi\leq k+1\}$. Taking then expectation, we have
\begin{multline}
\E \int_0^T\int_{\T^N}\Theta_k(u(x,t))|\gamma'(t)| dx dt+\E\int_{A_k^+} \gamma(t) dm(x,t,\xi)\\
= \frac{1}{2}\E\int_0^T\int_{\T^N} \gamma(t)\GG^2(x)\theta_k(u(x,t)) dx dt+\int_{\T^N}\Theta_k(u_0(x))dx.
\label{eq:kineticupreL1E}\end{multline}
Note that
\begin{equation}
(u-(k+1))^+\leq\Theta_k(u)\leq(u-k)^+,
\label{thetaplus}\end{equation}
for all $k\ge 0$, $u\in\R$. Note also
$$
\theta_{k_2}(u)\leq \mathbf{1}_{k_2\leq u}\leq\frac{(u-k_1)^+}{k_2-k_1},\quad 0\leq k_1<k_2.
$$
In particular, using that $m\geq 0$, \refe{eq:kineticupreL1E}, using \refe{D0} and taking $k_1=k$, $k_2=k_1+\theta-1$ where $\theta>1$ in \refe{eq:kineticupreL1E} gives
\begin{multline}
\E \int_0^T\int_{\T^N}(u(x,t)-(k+\theta))^+|\gamma'(t)| dx dt\\
\leq \frac{D_0}{2(\theta-1)}\E\int_0^T\int_{\T^N} (u(x,t)-k)^+\gamma(t)dx dt+\int_{\T^N} (u_0(x)-k)^+dx.
\label{stillgamma}\end{multline}
Choose $\theta$ large enough so that $ \frac{D_0}{2(\theta-1)}=\alpha<1$. Denote by $\psi_n(t)$ the function
$$
\psi_n(t)=\E \int_{\T^N}(u(x,t)-n\theta)^+dx
$$
and let $I_n\subset [0,T]$ be the set of Lebesgue points of $\psi_n$. Then $I=\cap_{n\in\N}I_n$ is of full measure. For $t\in I$, $t<T$, we take $k=n\theta$, $\gamma(s)=\min(1,\frac1\eps(s-(t+\eps))^-)$ for $\eps<T-t$ in \refe{stillgamma} and let $\eps\to0$. This yields the following inequality:
\begin{equation}\label{recpsin}
\psi_{n+1}(t)\leq\alpha\int_0^t\psi_n(s)ds+\psi_{n}(0).
\end{equation}
For $n=0$, using the bound 
$$ 
\psi_0(t)\leq\psi_1(t)+\theta\, T,
$$
\refe{recpsin} and Gronwall's Lemma, we obtain a bound
\begin{equation}\label{L1triv}
\psi_0(t)=\E\int_{\T^N} (u(x,t))^+ dx\leq M:=C(T,\|u_0\|_{L^1(\T^N)}),
\end{equation}
for $t\in I$. From \refe{recpsin} we also show recursively
$$
\psi_{n+1}(t)\leq\alpha^{n+1}\frac{t^{n+1}}{(n+1)!}M+\sum_{j=0}^n \alpha^j \frac{t^j}{j!}\psi_{n-j}(0).
$$
We deduce $\psi_{n+1}(t)\leq (M+1)e^T \delta(n+1)$, where 
$$
\delta(n)=\alpha^{n}+\sum_{j=0}^{n-1} \alpha^j \psi_{n-1-j}(0).
$$
Then $\lim_{n\to+\infty}\delta(n)=0$ since $\lim_{k\to+\infty}\psi_k(0)=0$. Set $\eps^+(k)=\delta([\frac{k}{\theta}])$, then  since the left hand side below is a decreasing function of $k$,
\begin{equation}\label{supEk}
\esssup_{t\in[0,T]}\E\int_{\T^N} (u(x,t)-k)^+ dx\leq\eps^+(k),
\end{equation}
and $\eps^+\colon\R_+\to\R_+$ is a function depending on $T$ and on the function $k\mapsto \|(u_0-k)^+\|_{L^1(\T^N)}$ only such that $\lim_{k\to+\infty}\eps^+(k)=0$. 
\medskip

{\bf Step 2.} The estimate on $m(A_k)$ follows from \refe{eq:kineticupreL1E}. Therefore, to conclude, we need to show an estimate on $\E(\esssup_{t\in[0,T]}\int_{\T^N} (u(x,t)-n\theta)^+ dx)$. This is the classical argument for semi-martingales that we will use. We will merely focus on the martingale term in \refe{eq:kineticupreL1gamma}. We first let $\gamma$ approach $\mathbf{1}_{(0,t)}$ as in Step 1. Then, by the Burkholder - Davis - Gundy inequality, we have, for a constant $C>0$, 
\begin{multline*}
\E\left[\esssup_{t\in[0,T]}\left|\sum_{j\geq 1}\int_0^t\int_{\T^N}g_j(x)\Theta_k'(u(x,t)) dx d\beta_k(t)\right|\right]\\
\leq C\E\left[\int_0^T\sum_{j\geq 1}\left(\int_{\T^N} |g_j(x)| |\Theta_k'(u(x,t))| dx\right)^2 \right]^{1/2} \\
\leq CD_0^{1/2} \E\left[\int_0^T\int_{\T^N} |\Theta_k'(u(x,t))|^2 dx \right]^{1/2},
\end{multline*}
where we have used the Cauchy-Schwarz inequality and \refe{D0} in the last line. Since
$$
|\Theta_k'(u)|^2\leq\mathbf{1}_{k\leq u}\leq (u-(k+1))^+
$$
and $\E X^{1/2}\leq (\E X)^{1/2}$ by Jensen's inequality, we obtain by \refe{supEk}
$$
\E\left[\esssup_{t\in[0,T]}\left|\sum_{j\geq 1}\int_0^t\int_{\T^N}g_j(x)\Theta_k'(u(x,t)) dx d\beta_k(t)\right|\right]\leq CD_0^{1/2}T\eps^+(k+1),
$$
as desired. This concludes the proof of the Proposition. \qed

\subsection{Uniqueness}\label{app:secuniqueness}

The condition at infinity {\it 2.} in Definition~\ref{def:kineticmeasureL1}, although weaker than the condition in \cite{DebusscheVovelle10} is enough to perform the proof of uniqueness, and of reduction of generalized solutions to solutions ({\it cf.} Theorem~15. in \cite{DebusscheVovelle10revised}). What is more problematic in the $L^1$ framework is the proof of pathwise continuity of solutions, {\it cf.} Corollary~16. in \cite{DebusscheVovelle10revised}. Indeed the proof of this result uses the equivalence between convergence to an equilibrium at the kinetic level and strong convergence at the level of functions. This is Lemma~7 in \cite{DebusscheVovelle10revised}, that we adapt here in the following way.

\begin{lemma}[Convergence to an equilibrium] Let $(X,\lambda)$ be a finite measure space. Let $(f_n)$ be a sequence of kinetic functions on $X\times\R$: $f_n(z,\xi)=\nu^n_z(\xi,+\infty)$ where $\nu^n$ are Young measures on $X$ satisfying
$$
\sup_n \int_X\int_\R |\xi|d\nu_z(\xi)d\lambda(z)<+\infty.
$$ 
Let $f$ be a kinetic function on $X\times\R$ such that $f_n\rightharpoonup f$ in $L^\infty(X\times\R)$ weak-*. Assume that $f_n$ and $f$ are equilibria: 
$$
f_n(z,\xi)=\mathbf{1}_{u_n(z)>\xi},\quad f(z,\xi)=\mathbf{1}_{u(z)>\xi}
$$
and assume that the following equi-integrability condition is satisfied:
\begin{equation}\label{equiintun}
\sup_n\|(u_n-k)^\pm\|_{L^1(X)}=o(1)\quad [k\to+\infty].
\end{equation}
Then $u_n\to u$ in $L^1(X)$ strong.
\label{lem:weakstrongeq}\end{lemma}

{\bf Proof.} We simply give the sketch of the proof of Lemma~\ref{lem:weakstrongeq}. Assume first 
\begin{equation}\label{uninfty}
R=\sup_n\|u_n\|_{L^\infty(X)}<+\infty.
\end{equation}
By testing the convergence $f_n\rightharpoonup f$ against $\varphi(z)\mathbf{1}_{\xi>-R}$, where $\varphi\in L^2(X)$, we obtain $u_n\rightharpoonup u$ in $L^2(X)$-weak. By testing the convergence $f_n\rightharpoonup f$ against $\xi\mathbf{1}_{\xi>-R}$, we obtain the convergence in norm $\|u_n\|_{L^2(X)}^2\to\|u\|_{L^2(X)}^2$ and since $L^2(X)$ is a Hilbert space, we conclude to $u_n\to u$ in $L^2(X)$ strong, and in particular in $L^1(X)$ strong. In the general case, we observe that the above arguments can be applied to $T_k(u_n)$, where $T_k$ is the truncation operator
\begin{equation}\label{truncate}
T_k(u)=\max(-k,\min(u,k)).
\end{equation}
We then use \refe{equiintun} to conclude to the convergence $u_n\to u$ in $L^1(X)$. \qed
\medskip

Proposition~\ref{prop:decaym} yields the equi-integrability estimate that we need. Therefore, as in \cite{DebusscheVovelle10revised}, we obtain the following two results.

\begin{theorem}[Uniqueness, Reduction] Let $u_0\in L^1(\T^N)$. Assume~\refe{D0}-\refe{D1}. 
Then, there is at most one solution in the sense of Definition~\ref{defkineticsolL1} with initial datum $u_0$ to \refe{stoSCL}. Besides, any generalized 
solution $f$ is actually a solution, {\it i.e.} if $f$ is a ge\-ne\-ra\-lized solution to \refe{stoSCL} with initial 
datum $\mathbf{1}_{u_0>\xi}$, then there exists a solution $u$ to \refe{stoSCL} with initial datum 
$u_0$ such that $f(x,t,\xi)=\mathbf{1}_{u(x,t)>\xi}$ a.s., for a.e. $(x,t,\xi)$. 
\label{th:UaddL1}\end{theorem}

\begin{theorem}[Continuity in time]Let $u_0\in L^1(\T^N)$. Assume~\refe{D0}-\refe{D1} are satisfied. Then the solution $u$ to \refe{stoSCL} with initial datum $u_0$ has a representative in the space $L^1(\Omega;L^\infty(0,T;L^1(\T^N)))$ with almost sure continuous trajectories in $L^1(\T^N)$.
\label{cor:timecontinuityL1}\end{theorem}

\subsection{Existence}

\begin{theorem}[Existence] Let $u_0\in L^1(\T^N)$. Assume~\refe{D0}-\refe{D1}. Then, there exists a solution in the sense of Definition~\ref{defkineticsolL1} with initial datum $u_0$ to \refe{stoSCL}. 
\label{th:existL1}\end{theorem}

{\bf Proof.} We approach $u_0$ by $u^n_0:=T_n(u_0)$, where the truncation operator is defined in \refe{truncate}. By \cite{DebusscheVovelle10} this defines a sequence of solutions $(u_n)$ and by the
contraction property in $L^1$ it is a Cauchy sequence, hence converges to a $u\in L^1(\Omega\times\T^N\times(0,T))$. This $u$ is predictable. Let now $m^n$ be the kinetic measure associated to $u^n$. By Proposition~\ref{prop:decaym}, we have 
\begin{equation}
\sup_n \E m^n(K_r)\leq C_r
\label{estimmn}\end{equation}
for $r\in\N^*$, where $K_r:=\T^N\times[0,T]\times[-r,r]$. Let $\mathcal{M}_r$ denote the space of bounded Borel measures over $K_r$ (with norm given by the total variation of measures). It is the topological dual of $C(K_r)$, the set of 
continuous functions on $K_r$. Since $\mathcal{M}_r$ is separable ($C(K_r)$ is) the space $L^1(\Omega;\mathcal{M}_r)$ is 
the topological dual space of $L^1(\Omega,C(K_r))$, {\it c.f.} Th\'eor\`eme~1.4.1 in \cite{Droniou01}. The estimate \refe{estimmn} gives a uniform bound on $(m^n)$ in $L^1(\Omega,\mathcal{M}_r)$: there exists $m_r\in L^1(\Omega,\mathcal{M}_r)$ such that  
up to subsequence, $m^n\rightharpoonup m_r$ in $L^1(\Omega;\mathcal{M}_r)$-weak star. By a dia\-go\-nal process, we obtain, for $r\in\N^*$, $m_r=m_{r+1}$ in $L^1(\Omega;\mathcal{M}_r)$ and the convergence in all the spaces $L^1(\Omega;\mathcal{M}_r)$-weak star of a single subsequence still denoted $(m^n)$. Let us then set $m=m_r$ on $K_r$, a.s. The conditions {\it 1.} and {\it 3.} in Definition~\ref{def:kineticmeasureL1} are stable by weak convergence, hence satisfied by $m$. We deduce that condition~{\it 2.} is satisfied thanks to the uniform estimate of Proposition~\ref{prop:decaym}. This shows that $u$ is a solution to \refe{stoSCL} with initial datum $u_0$. \qed
\medskip

We may also prove the existence of solution by convergence of the parabolic approximation as in \cite{DebusscheVovelle10}. We will not give the details of the proof. Our final result is the following one.

\begin{theorem}[Resolution of \refe{stoSCL} in $L^1$] Let $u_0\in L^1(\T^N)$. There exists a unique measurable $u\colon\T^N\times [0,T]\times\Omega\to\R$ solution to~\refe{stoSCL} with initial datum $u_0$ in the sense of Definition~\ref{defkineticsolL1}. Besides, $u$ has almost surely continuous trajectories in $L^1(\T^N)$ and $u$ is the a.s. limit in $L^1(\T^N\times(0,T))$ of the parabolic approximation to \refe{stoSCL}.  Moreover, given $u_0^1$ and $u_0^2\in L^1(\T^N)$, the following holds:
$$
\|u^1(t)-u^2(t)\|_{L^1(\T^N)}\le \|u_0^1-u_0^2\|_{L^1(\T^N)},\; a.s.
$$
\end{theorem}

\begin{remark}[Multiplicative noise] Here we have developed the $L^1$ theory for an additive noise (the functions $g_k$ are independent on $u$) since this was assumed from the start, but all the statement above remain true if the noise is multiplicative with the same hypotheses as in \cite{DebusscheVovelle10} for example.
\end{remark}

\end{appendix}

\providecommand{\bysame}{\leavevmode\hbox to3em{\hrulefill}\thinspace}
\providecommand{\MR}{\relax\ifhmode\unskip\space\fi MR }
\providecommand{\MRhref}[2]{%
  \href{http://www.ams.org/mathscinet-getitem?mr=#1}{#2}
}
\providecommand{\href}[2]{#2}

\end{document}